\documentclass[11pt]{article}

\usepackage{amsmath}
\usepackage{amsthm}
\usepackage{amssymb}
\usepackage{amscd}

\usepackage{a4wide}

\theoremstyle{definition}
\newtheorem{theorem}{Theorem}[section]
\newtheorem{prop}[theorem]{Proposition}
\newtheorem{lemma}[theorem]{Lemma}
\newtheorem{corollary}[theorem]{Corollary}
\newtheorem{definition}[theorem]{Definition}

\newtheorem{remark}[theorem]{Remark}
\newtheorem{conjecture}[theorem]{Conjecture}

\renewcommand\hat{\widehat}
\renewcommand\tilde{\widetilde}
\newcommand\ep{\varepsilon}
\newcommand\KK{\mathcal{K}}
\newcommand\JJ{\mathcal{J}}

\newcommand\PP{\mathcal{P}}

\newcommand\Comp{\mathbb{C}}
\newcommand\Real{\mathbb{R}}
\newcommand\Rat{\mathbb{Q}}
\newcommand\Int{\mathbb{Z}}

\newcommand\rank{\operatorname{rk}}

\newenvironment{demo}[1]{%
  \trivlist
  \item[\hskip\labelsep
        {\bf #1.}]
}{%
  \endtrivlist
}

\allowdisplaybreaks

\title{
Birational Rowmotion and Coxeter-motion \\
on Minuscule Posets
}

\author{
Soichi OKADA
\footnote{
Graduate School of Mathematics, Nagoya University, 
Furo-cho, Chikusa-ku, Nagoya 464-8602, Japan, 
{\tt okada@math.nagoya-u.ac.jp}
}
}

\date{
}

\begin{document}

\maketitle

\begin{abstract}
Birational rowmotion is a discrete dynamical system 
on the set of all positive real-valued functions on a finite poset, 
which is a birational lift of combinatorial rowmotion on order ideals.
It is known that combinatorial rowmotion for a minuscule poset has order equal 
to the Coxeter number, 
and exhibits the file homomesy phenomenon for refined order ideal cardinality statistic.
In this paper we generalize these results to the birational setting.
Moreover, as a generalization of birational promotion on a product of two chains, 
we introduce birational Coxeter-motion on minuscule posets, and prove that it enjoys 
periodicity and file homomesy.
\par
Mathematics Subject Classification (MSC2010): 05E18 (primary), 06A11 (secondary)
\end{abstract}

\section{%
Introduction
}

Rowmotion (at the combinatorial level) is a bijection $R$ on the set $\JJ(P)$ 
of order ideals of a finite poset $P$, which assigns to $I \in \JJ(P)$ 
the order ideal $R(I)$ generated by the minimal elements of the complement $P \setminus I$.
The map $R$ can be also described in terms of toggles.
For each $v \in P$, let $t_v : \JJ(P) \to \JJ(P)$ be the map given by
\begin{equation}
\label{eq:Ctoggle}
t_v(I)
 =
\begin{cases}
I \cup \{ v \} &\text{if $v \not\in I$ and $I \cup \{ v \} \in \JJ(P)$,} \\
I \setminus \{ v \} &\text{if $v \in I$ and $I \setminus \{ v \} \in \JJ(P)$,} \\
I &\text{otherwise,}
\end{cases}
\end{equation}
and call it the \emph{toggle} at $v$, 
Then the rowmotion map $R$ is expressed as the composition
\begin{equation}
\label{eq:row=toggle}
R = t_{v_1} \circ t_{v_2} \circ \cdots \circ t_{v_N},
\end{equation}
where $(v_1, v_2, \dots, v_N)$ is any linear extension of $P$, 
i.e., a list of elements of $P$ such that $v_i < v_j$ in $P$ implies $i<j$.
This rowmotion has been studied from several perspectives and under various names.
See \cite{SW} and \cite{TW} for the history and references.

Rowmotion exhibits nice properties such as periodicity and homomesy 
on special posets including root posets (see \cite{Pan, AST}) 
and minuscule posets (see \cite{RS, RW}).
In general, given a set $S$ and a bijection $f : S \to S$, 
we say that a statistic $\theta : S \to \Real$ is \emph{homomesic} with respect to $f$ 
if there exists a constant $C$ such that for any $\langle f \rangle$-orbit $T$
$$
\frac{ 1 }{ \# T }
\sum_{x \in T} \theta(x)
 =
C.
$$
We refer the reader to \cite{R} for the homomesy phenomenon.
For a minuscule poset $P$ and a simple root $\alpha \in \Pi$, we put
\begin{equation}
\label{eq:file}
P^\alpha = \{ v \in P : c(v) = \alpha \},
\end{equation}
where $c : P \to \Pi$ is the coloring of $P$ with color set $\Pi$, the set of simple roots. 
This subset $P^\alpha$ is called the \emph{file} corresponding to $\alpha$.
(See Section~3 for a definition of minuscule posets and related terminologies.)

If $P$ is a minuscule poset, then the associated rowmotion map $R$ has the following properties:

\begin{theorem}
\label{thm:prototype}
Let $P$ be a minuscule poset associated to a minuscule weight $\lambda$ of a simple Lie algebra $\mathfrak{g}$.
Then we have
\begin{enumerate}
\item[(a)]
(periodicity, Rush--Shi \cite[Thoerem~1.4]{RS})
The rowmotion map $R$ has a finite order equal to the Coxeter number $h$ of $\mathfrak{g}$.
\item[(b)]
(file homomesy, Rush--Wang \cite[Theorem~1.2]{RW})
For each simple root $\alpha \in \Pi$, 
the refined order ideal cardinality $\# (I \cap P^\alpha)$ is homomesic with respect to $R$.
More precisely, for any $I \in \JJ(P)$, we have
$$
\frac{1}{h}
 \sum_{k=0}^{h-1} \# \left( R^k(I) \cap P^\alpha \right)
 =
\langle \varpi^\vee, \lambda \rangle,
$$
where $\varpi^\vee$ is the fundamental coweight corresponding to $\alpha$.
\end{enumerate}
\end{theorem}

One motivation of this paper is to lift the results in the above theorem to the birational level.

Einstein--Propp \cite{EP} introduced birational rowmotion 
by lifting the notion of toggles from the combinatorial level to the piecewise-linear level, 
and then to the birational level. 
Given a finite poset $P$, let $\hat{P} = P \sqcup \{ \hat{1}, \hat{0} \}$ be the poset obtained from $P$ 
by adjoining an extra maximum element $\hat{1}$ and an extra minimum element $\hat{0}$.
For positive real numbers $A$ and $B$, we put
$$
\KK^{A,B}(P)
 =
\{ F : \hat{P} \to \Real_{>0} \mid F(\hat{1}) = A, \ F(\hat{0}) = B \},
$$
where $\Real_{>0}$ denotes the set of positive real numbers.
For $v \in P$, we define the \emph{birational toggle} $\tau^{A,B}_v : \KK^{A,B}(P) \to \KK^{A,B}(P)$ at $v$ by
\begin{equation}
\label{eq:Btoggle}
\left( \tau^{A,B}_v F \right)(x)
 =
\begin{cases}
\dfrac{ 1 }
      { F(v) }
\cdot
\dfrac{ \sum_{w \in \hat{P}, \, w \lessdot v} F(w) }
      { \sum_{z \in \hat{P}, \, z \gtrdot v} 1/F(z) }
 &\text{if $x = v$,}
\\
F(x) &\text{otherwise,}
\end{cases}
\end{equation}
where the symbol $x \gtrdot y$ means that $x$ covers $y$, 
i.e., $x>y$ and there is no element $z$ such that $x>z>y$.
It is clear that $\tau^{A,B}_v$ is an involution.
(See (\ref{eq:PLtoggle}) for a definition of piecewise-linear toggles.)
Then we define \emph{birational rowmotion} $\rho^{A,B} : \KK^{A,B}(P) \to \KK^{A,B}(P)$ by
\begin{equation}
\label{eq:Brow}
\rho^{A,B}
 = 
\tau^{A,B}_{v_1} \circ \cdots \circ \tau^{A,B}_{v_N},
\end{equation}
where $(v_1, \dots, v_N)$ is a linear extension of $P$.
It can be shown that the definition of $\rho^{A,B}$ is independent of the choices of linear extensions.
Since rowmotion is defined by toggling from top to bottom, 
we have a recursion formula for the values of the birational rowmotion map:
\begin{equation}
\label{eq:Brow_inductive}
\left( \rho^{A,B} F \right)(v)
 =
\frac{ 1 }
     { F(v) }
\cdot
\frac{ \sum_{w \in \hat{P}, \, w \lessdot v} F(w) }
     { \sum_{z \in \hat{P}, \, z \gtrdot v} 1/\left( \rho^{A,B} F \right) (z) }.
\end{equation}
We omit the superscript ${}^{A,B}$ and simply write $\KK(P)$, $\tau_v$ and $\rho$ 
when there is no confusion.

For birational rowmotion on a product of two chains, 
periodicity and (multiplicative) file homomesy are obtained by 
Grinberg--Roby \cite{GR2} and Einstein--Propp \cite{EP} respectively.
In this paper we generalize their results from products of two chains (type $A$ minuscule posets) 
to arbitrary minuscule posets.

For a minuscule poset and a simple root $\alpha \in \Pi$, we define
\begin{equation}
\label{eq:Phi}
\Phi_\alpha(F)
 = 
\prod_{v \in P^\alpha} F(v)
\end{equation}
for $F \in \KK^{A,B}(P)$.
Our main results for birational rowmotion are summarized as follows:

\begin{theorem}
\label{thm:main1}
Let $P$ be the minuscule poset associated to a minuscule weight $\lambda$ 
of a finite dimensional simple Lie algebra $\mathfrak{g}$.
Let $\rho = \rho^{A,B}$ be the birational rowmotion map.
Then we have
\begin{enumerate}
\item[(a)]
(periodicity)
The map $\rho$ 
has finite order equal to the Coxeter number $h$ of $\mathfrak{g}$.
\item[(b)]
(reciprocity)
For any $v \in P$ and $F \in \KK^{A,B}(P)$, we have
\begin{equation}
\label{eq:reciprocity}
\left( \rho^{\rank(v)} F \right) (v)
 =
\frac{ AB }
     { F( \iota v) },
\end{equation}
where $\rank : P \to \{ 1, 2, \dots,  h-1 \}$ is the rank function of the graded poset $P$ 
and $\iota : P \to P$ is the canonical involutive anti-automorphism of $P$ 
(see Proposition~\ref{prop:involution}).
\item[(c)]
(file homomesy)
For a simple root $\alpha$, we have
\begin{equation}
\label{eq:homomesyR}
\prod_{k=0}^{h-1} \Phi_\alpha (\rho^k F)
 =
A^{h \langle \varpi^\vee, -w_0 \lambda \rangle}
B^{h \langle \varpi^\vee, \lambda \rangle}
\end{equation}
for any $F \in \KK^{A,B}(P)$, 
where $w_0$ is the longest element of the Weyl group $W$ of $\mathfrak{g}$, 
and $\varpi^\vee$ is the fundamental coweight corresponding to $\alpha$.
\end{enumerate}
\end{theorem}

Part (a) of this theorem is established in \cite{GR1,GR2} 
except for the type $E_7$ minuscule poset.
In this paper we provide a way to settle the $E_7$ case by using a computer.
For a type $A$ minuscule poset,
Part (b) is obtained in \cite[Theorem~30]{GR2}.
Our proof of Part (b) is based on a case-by-case analysis (with a help of computer in types $E_6$ and $E_7$).
Part (c) in type $A$ follows from Einstein--Propp \cite[Theorems~7.3 and 8.5]{EP} 
(see Musiker--Roby \cite[Theorem~2.16]{MR} for another proof).
We will give an almost uniform proof to Part (c).
Also we can use tropicalization (or ultradiscretization) to deduce the results 
for combinatorial rowmotion in Theorem~\ref{thm:prototype} (see Section~2).

Another aim of this paper is to introduce and study birational Coxeter-motion on minuscule posets, 
which is regarded as a generalization of birational promotion on a product of two chains 
(see \cite[Definition~5.3]{EP}).
For a simple root $\alpha \in \Pi$, we define $\sigma^{A,B}_\alpha : \KK^{A,B}(P) \to \KK^{A,B}(P)$ 
as the composition
\begin{equation}
\label{eq:sigma}
\sigma^{A,B}_\alpha = \prod_{v \in P_\alpha} \tau^{A,B}_v, 
\end{equation}
which is independent of the order of composition.
Then a \emph{Coxeter-motion map} is a product of all the $\sigma^{A,B}_\alpha$'s in any order.
Our results for birational Coxeter-motion are stated as follows:

\begin{theorem}
\label{thm:main2}
Let $P$ be the minuscule poset.
Let $\gamma = \gamma^{A,B}$ be a birational Coxeter-motion map.
Then we have
\begin{enumerate}
\item[(a)]
(periodicity)
The map $\gamma$ has finite order equal to the Coxeter number $h$.
\item[(b)]
(file homomesy)
For each simple root $\alpha \in \Pi$, we have
\begin{equation}
\prod_{k=0}^{h-1} \Phi^\alpha (\gamma^k F)
 =
A^{h \langle \varpi^\vee, -w_0 \lambda \rangle}
B^{h \langle \varpi^\vee, \lambda \rangle}.
\end{equation}
\end{enumerate}
\end{theorem}

If $P$ is a type $A$ minuscule poset and $\pi$ is the birational promotion map 
(a special case of birational Coxeter-motion maps), 
then there is an explicitly defined ``recombination map'' $\mathfrak{R}$ 
such that $\mathfrak{R} \rho = \pi \mathfrak{R}$ (see \cite[Theorem~8.2]{EP}), 
which, together with Theorem~\ref{thm:main1} (a), implies Part (a) of the above theorem. 
We prove Part (a) for arbitrary minuscule posets by showing that any birational Coxeter-motion map 
is conjugate to the birational rowmotion map in the birational toggle group (Theorem~\ref{thm:conj} below). 
By applying tropicalization to Part (a), 
we obtain the periodicity of piecewise-linear promotion, 
which is proved in \cite[Theorem~1.12]{GPT} via quiver representation.
Part (b) in type $A$ is obtained in \cite[Theorem~7.3]{EP}.

Hopkins \cite{Ho} obtains another example of homomesy for the birational rowmotion 
for a wider class of posets including minuscule posets.

\begin{theorem}
\label{thm:Hopkins}
(Hopkins \cite[Theorem~4.43]{Ho})
Let $P$ be a minuscule poset and $\rho = \rho^{A,B}$ the birational rowmotion map.
For $F \in \KK^{A,B}(P)$, we define
$$
\Psi(F)
 = 
\prod_{x \in P}
 \frac{ F(x) }
      { \sum_{y \in \hat{P}, y \lessdot x} F(y) }.
$$
Then we have
$$
\prod_{k=0}^{h-1} \Psi(\rho^k F) 
 =
\left( \frac{ A }{ B} \right)^{\# P}.
$$
\end{theorem}

Via tropicalization, this theorem reduces to the homomesy phenomenon 
of the antichain cardinality statistic, which was proved in \cite[Theorem~1.4]{RW}.
In a forthcoming paper \cite{O}, we use explicit formulas for iterations of the birational rowmotion map 
to give refinements of Theorem~\ref{thm:Hopkins}.
Our refinement in type $A$ provides a birational lift of the homomesy given in \cite[Proof of Theorem~27]{EP}.

The remaining of this paper is organized as follows.
We collect some general facts concerning birational rowmotion in Section~2, 
and give a definition and properties of minuscule posets in Section~3.
In Sections 4 to 6 we give a proof of our main theorems.
The periodicity in Theorem~\ref{thm:main1} (a) and Theorem~\ref{thm:main2} (a) is proved in Section~4, 
and the reciprocity in Theorem~\ref{thm:main1} (b) is verified in Section~5.
In Section~6, after investigating local properties around a file, 
we complete the proof of file homomesy in Theorem~\ref{thm:main1} (c) and Theorem~\ref{thm:main2} (b).

\subsection*{%
Acknowledgements
}
This work was partially supported by 
JSPS Grants-in-Aid for Scientific Research No.~18K03208. 
The author is grateful to Tom Roby for fruitful discussions.
\section{%
Generalities on rowmotion
}

In this section, we explain how combinatorial and birational rowmotion are related, 
and give some general facts about birational rowmotion.

\subsection{%
Combinatorial, piecewise-linear and birational rowmotion
}

We begin with recalling the definition of piecewise-linear toggles and rowmotion.
Given a finite poset $P$ and real numbers $a$, $b$, we put
$$
\PP^{a,b}(P)
 =
\{ f : \hat{P} \to \Real : f(\hat{1}) = a, \, f(\hat{0}) = b \},
$$
where $\hat{P} = P \sqcup \{ \hat{1}, \hat{0} \}$.
We define the \emph{piecewise-linear toggles} 
$\tilde{t}^{\pm,a,b}_v : \PP^{a,b}(P) \to \PP^{a,b}(P)$ at $v \in P$ 
by the formulas
\begin{equation}
\label{eq:PLtoggle}
\begin{aligned}
\left( \tilde{t}^{+,a,b}_v f \right)(v)
 &=
\max \{ f(w) : w \in \hat{P}, \, w \lessdot v \} + \min \{ f(z) : z \in \hat{P}, \, z \gtrdot v \} - f(v),
\\
\left( \tilde{t}^{-,a,b}_v f \right)(v)
 &=
\min \{ f(w) : w \in \hat{P}, \, w \lessdot v \} + \max \{ f(z) : z \in \hat{P}, \, z \gtrdot v \} - f(v),
\end{aligned}
\end{equation}
and $\left( \tilde{t}^{\pm,a,b}_v f \right)(x) = f(x)$ for $x \neq v$.
For an order ideal $I \in \JJ(P)$, let $\chi^\pm_I$ be the characteristic functions defined by
$$
\chi^+_I(v)
 = 
\begin{cases}
 0 &\text{if $v \in I$ or $v = \hat{0}$,} \\
 1 &\text{if $v \in P \setminus I$ or $v = \hat{1}$,}
\end{cases}
\quad
\chi^-_I(v)
 = 
\begin{cases}
 1 &\text{if $v \in I$ or $v = \hat{0}$,} \\
 0 &\text{if $v \in P \setminus I$ or $v = \hat{1}$.}
\end{cases}
$$
Then it follows from the definition (\ref{eq:Ctoggle}) and (\ref{eq:PLtoggle}) that 
the toggle $\tilde{t}^{\pm,a,b}_v$ is a piecewise-linear lift of the combinatorial toggle $t_v$ 
in the following sense:
\begin{equation}
\label{eq:C-PL-toggle}
\tilde{t}^{+,1,0}_v ( \chi^+_I )
 =
\chi^+_{t_v I},
\quad
\tilde{t}^{-,0,1}_v ( \chi^-_I )
 =
\chi^-_{t_v I}.
\end{equation}
The \emph{piecewise-linear rowmotion} map $\tilde{R}^{\pm,a,b} : \PP^{a,b}(P) \to \PP^{a,b}(P)$ 
is defined by
$$
\tilde{R}^{\pm,a,b} = \tilde{t}^{\pm,a,b}_{v_1} \circ \cdots \circ \tilde{t}^{\pm,a,b}_{v_N},
$$
where $(v_1, \dots, v_N)$ is a linear extension of $P$.

A rational function $F(X_1, \cdots, X_m) \in \Rat(X_1, \cdots, X_m)$ is called 
\emph{subtraction-free} if 
$F$ is expressed as a ratio $F = G/H$ of two polynomials $G(X_1, \cdots, X_m)$ and $H(X_1, \cdots, X_m) 
\in \Int[X_1, \dots, X_m]$ with nonnegative integer coefficients.
By using
$$
\lim_{\ep \to +0}
 \ep \log (e^{a/\ep} + e^{b/\ep})
 =
\max \{ a, b \},
\quad
\lim_{\ep \to -0}
 \ep \log (e^{a/\ep} + e^{b/\ep})
 =
\min \{ a, b \},
$$
we can see that, if $F(X_1, \dots, X_m)$ is subtraction-free, 
then for any real numbers $x_1, \dots, x_m \in \Real$ 
the limits
$$
f^\pm(x_1, \cdots, x_m)
 = 
\lim_{\ep \to \pm 0}
 \ep \log F(e^{x_1/\ep}, \cdots, e^{x_m/\ep})
$$
exist and 
$f^+(x_1, \dots, x_m)$ (resp. $f^-(x_1, \dots, x_m)$) is the piecewise-linear function in $x_1, \dots, x_m$ 
obtained from $F$ by replacing the multiplication $\cdot$, the division $/$ and the addition $+$ 
with the addition $+$, the subtraction $-$ and the maximum $\max$ (resp. the minimum $\min$).
This procedure from $F$ to $f^\pm$ are called the tropicalization (or ultradiscretization).

\begin{prop}
\label{prop:tropical}
Let $P$ be a finite poset.
Let $R : \JJ(P) \to \JJ(P)$ and $\rho = \rho^{A,B} : \KK^{A,B}(P) \to \KK^{A,B}(P)$ 
be the combinatorial and birational rowmotion maps respectively.
Let $m : P \times \Int \to \Int$ be a map with finite support.
If there is a integers $p$ and $q$ such that
\begin{equation}
\label{eq:Btropical}
\prod_{(v,k) \in P \times \Int}
 \left[ \left( \rho^k F \right)(v) \right]^{m(v,k)}
 =
A^p B^q
\end{equation}
for any $F \in \KK^{A,B}(P)$, then
\begin{equation}
\label{eq:Ctropical}
\sum_{(v,k) \in P \times \Int} 
 m(v,k) \chi[v \not\in R^k(I)] 
 =
p,
\quad
\sum_{(v,k) \in P \times \Int} 
 m(v,k) \chi[v \in R^k(I)] 
 =
q,
\end{equation}
where $\chi[S] = 1$ if $S$ is true and $0$ if $S$ is false.
\end{prop}

\begin{demo}{Proof}
By applying the tropicalization procedure to (\ref{eq:Btropical}), we obtain
$$
\sum_{(v,k) \in P \times \Int} 
 m(v,k) \left( \tilde{R}^{\pm,a,b} f \right)(v) 
 =
a p + b q
$$
for any $f \in \PP^{a,b}(P)$.
Then specializing $f = \chi^\pm_I$ and using (\ref{eq:C-PL-toggle}), we obtain (\ref{eq:Ctropical}).
\qed
\end{demo}

\begin{corollary}
\label{cor:tropical}
\begin{enumerate}
\item[(a)]
If $\left( \rho^h F \right)(v) = F(v)$ for any $F \in \KK^{A,B}(P)$ and $v \in P$, 
then $R^h (I) = I$ any $I \in \JJ(P)$.
\item[(b)]
Let $v$ and $w \in P$ and $k$ a positive integer.
If $\left( \rho^k F \right)(v) \cdot F(w) = AB$, 
then $v \in R^k(I)$ and $w \not\in I$ are equivalent for any $I \in \JJ(P)$.
\item[(c)]
Let $M$ be a subset of $P$ and $h$ be a positive integer.
If $\prod_{k=0}^{h-1} \prod_{v \in M} \left( \rho^k F \right)(v) = A^p B^q$ 
for any $F \in \KK^{A,B}(P)$, 
then we have $\sum_{k=0}^{h-1} \# \left( R^k(I) \cap M \right) = q$ for any $I \in \JJ(P)$.
\end{enumerate}
\end{corollary}

Similar statements hold for birational Coxter-motion.

\subsection{%
Birational rowmotion on graded posets
}

In this subsection we present some properties of birational rowmotion on graded posets.
A poset $P$ is called \emph{graded of height $n$} 
if there exists a rank function $\rank : P \to \{ 1, 2, \dots, n \}$ satisfying the following three conditions:
\begin{enumerate}
\item[(i)]
If $v$ is minimal in $P$, then $\rank(v) = 1$;
\item[(ii)]
If $v$ is maximal in $P$, then $\rank(v) = n$;
\item[(iii)]
If $v$ covers $w$, then $\rank(v) = \rank(w)+1$.
\end{enumerate}

\begin{lemma}
\label{lem:order}
If $P$ is a graded poset of height $n$ 
and the birational rowmotion map $\rho^{A,B}$ has a finite order $N$, 
then $N$ is divisible by $n+1$.
\end{lemma}

\begin{demo}{Proof}
By Corollary~\ref{cor:tropical} (a), 
we have $R^N(I) = I$ for all $I \in \JJ(P)$.
On the other hand, it is easy to see that the $\langle R \rangle$-orbit of the empty order ideal $\emptyset$ 
has length $n+1$.
Hence we see that $n+1$ divides $N$.
\qed
\end{demo}

The following lemma gives a relation between $\rho^{A,B}$ and $\rho^{1,1}$.

\begin{lemma}
\label{lem:A=B=1}
Let $P$ be a graded poset of height $n$.
For a map $F : P \to \Real_{>0}$ and positive real numbers $A$, $B \in \Real_{>0}$, 
we denote by $F^{A,B} \in \KK^{A,B}(P)$ the extension of $F$ to $\hat{P}$ such that 
$F^{A,B}(\hat{1}) = A$ and $F^{A,B}(\hat{0}) = B$.
For $1 \le k \le n$ and $v \in P$, we have
\begin{equation}
\left(
 \left( \rho^{A,B} \right)^k F^{A,B}
\right)(v)
 =
\left(
 \left( \rho^{1,1} \right)^k F^{1,1}
\right)(v)
\times
\begin{cases}
 A &\text{if $1 \le k \le \rank(v)-1$,} \\
 AB &\text{if $k = \rank(v)$,} \\
 B &\text{if $\rank(v)+1 \le k \le n$,} \\
 1 &\text{if $k = n+1$.}
\end{cases}
\end{equation}
\end{lemma}

\begin{demo}{Proof}
We can use the recursive formula (\ref{eq:Brow_inductive}) 
to proceed by the double induction on $k$ and $n-\rank(v)$.
\qed
\end{demo}

\subsection{%
Change of variables
}

Let $P$ be a finite poset.
Given an initial state $X \in \KK^{A,B}(P)$, 
we regard $X(v)$ ($v \in P$) as indeterminates.
In the computation of $\left( \rho^k X \right)(v)$ ($v \in P$) 
of iterations of the birational rowmotion map $\rho = \rho^{A,B}$, 
it is convenient to change variables from $\{ X(v) : v \in P \}$ 
to $\{ Z(v) : v \in P \}$ defined by the formula
\begin{equation}
\label{eq:X2Z}
Z(v)
 = 
\begin{cases}
 X(v) &\text{if $v$ is minimal,} \\
 \dfrac{ X(v) }
       { \sum_{w \in P, \, w \lessdot v} X(w) }
 &\text{otherwise.}
\end{cases}
\end{equation}
This change of variables is used in \cite{MR} to describe a lattice path formula 
for birational rowmotion on a type $A$ minuscule poset.
Then the inverse change of variables is given by
\begin{equation}
\label{eq:Z2X}
X(v) = \sum Z(v_1) Z(v_2) \cdots Z(v_r),
\end{equation}
where the sum is taken over all saturated chains $v_1 \gtrdot \cdots \gtrdot v_r$ in $P$ 
such that $v_1 = v$ and $v_r$ is minimal in $P$.
Note that this change of variables is a birational lift of Stanley's transfer map 
between the order polytope and the chain polytope of a poset (see \cite[Section~3]{Sta}).

\section{
Minuscule posets
}

In this section we review a definition and properties of minuscule posets.

\subsection{%
Definition and properties of minuscule posets
}

Let $\mathfrak{g}$ be a finite dimensional simple Lie algebra over the complex number field $\Comp$ 
of type $X_n$, where $X \in \{ A, B, C, D, E, F, G \}$ and $n$ is the rank of $\mathfrak{g}$.
We fix a Cartan subalgebra $\mathfrak{h}$ and choose a positive root system $\Delta_+$ 
of the root system $\Delta \subset \mathfrak{h}^*$.
Let $\Pi = \{ \alpha_1, \dots, \alpha_n \}$ be the set of simple roots, 
where we follow \cite[Planche~I--IX]{B1} for the numbering of simple roots.
We denote by $\varpi_i$ the fundamental weight corresponding to 
the $i$th simple root $\alpha_i$.
Let $\Delta^\vee_+ \subset \mathfrak{h}$ be the positive coroot system.
Let $W$ be the Weyl group of $\mathfrak{g}$, which acts on $\mathfrak{h}$ and $\mathfrak{h}^*$.
The simple reflections $\{ s_\alpha : \alpha \in \Pi \}$ generate $W$.

For a dominant integral weight $\lambda$, 
we denote by $V_{X_n,\lambda}$ the irreducible $\mathfrak{g}$-module with highest weight $\lambda$ 
and by $L_{X_n,\lambda}$ the set of weights of $V_{X_n,\lambda}$.
We say that $\lambda$ is \emph{minuscule} if $L_{X_n,\lambda}$ is a single $W$-orbit.
See \cite[VIII, \S7, n${}^\circ$3]{B2} for properties of minuscule weights.
It is known that minuscule weights are fundamental weights.
Table~\ref{tab:minuscule} is the list of minuscule weights.
\begin{table}[ht]
\caption{List of minuscule weights}
\label{tab:minuscule}
\centering
\begin{tabular}{c|c|c}
type & minuscule weights & Coxeter number \\
\hline
$A_n$ & $\varpi_1, \varpi_2, \dots, \varpi_n$ & $n+1$ \\
$B_n$ & $\varpi_n$ & $2n$ \\
$C_n$ & $\varpi_1$ & $2n$ \\
$D_n$ & $\varpi_1, \varpi_{n-1}, \varpi_n$ & $2n-2$ \\
$E_6$ & $\varpi_1, \varpi_6$ & $12$ \\
$E_7$ & $\varpi_7$ & $18$ \\
$E_8$ & none & $30$ \\
$F_4$ & none & $12$ \\
$G_2$ & none & $6$
\end{tabular}
\end{table}

Let $\lambda$ be a minuscule weight of a simple Lie algebra $\mathfrak{g}$ of type $X_n$.
We equip the set of wegiths $L_{X_n,\lambda}$ with a poset structure by defining 
$\mu \ge \nu$ if $\nu - \mu$ is a linear combination of simple roots 
with nonnegative integer coefficients.
We note that $\lambda$ is the minimum element of the poset $L_{X_n,\lambda}$.

\begin{definition}
Let $\mathfrak{g}$ be a simple Lie algebra of type $X_n$ 
and $\lambda$ a minuscule weight.
Then the \emph{minuscule poset} $P_{X_n,\lambda}$ is defined by
\begin{equation}
P_{X_n,\lambda}
 = 
\{ \beta^\vee \in \Delta^\vee_+ : \langle \beta^\vee, \lambda \rangle = 1 \},
\end{equation}
where the partial ordering on $P_{X_n,\lambda}$ is given by saying that 
$\alpha^\vee \ge \beta^\vee$ if $\alpha^\vee - \beta^\vee$ 
is a linear combination of simple coroots with nonnegative integer coefficients.
\end{definition}

\begin{prop}
\label{prop:minuscule}
Let $\lambda$ be a minuscule weight 
and $P_{X_n,\lambda}$ be the corresponding minuscule poset.
Then we have
\begin{enumerate}
\item[(a)]
(\cite[Propotisions~3.2, 4.1]{P})
The poset $L_{X_n,\lambda}$ is a distributive lattice.
\item[(b)]
(\cite[Theorem~11]{P})
There exists a unique map $c : P_{X_n,\lambda} \to \Pi$, called the \emph{coloring} of $P_{X_n,\lambda}$, 
such that the map
$$
\JJ(P_{X_n,\lambda}) \ni I \mapsto \lambda - \sum_{v \in I} c(v) \in L_{X_n,\lambda}
$$
gives an isomorphism of posets.
\end{enumerate}
\end{prop}

If $\lambda$ is a minuscule weight, then the stabilizer $W_\lambda$ of $\lambda$ in $W$ 
is the maximal parabolic subgroup generated by $\{ s_\beta : \beta \in \Pi \setminus \{ \alpha \} \}$, 
where $\alpha$ is the simple root corresponding to the fundamental weight $\lambda$.

\begin{prop}
\label{prop:involution}
Let $P_{X_n,\lambda}$ be the minuscule poset corresponding to a minuscule weight $\lambda$, 
and $w_\lambda$ the longest element of the stabilizer $W_\lambda$.
Then the map
$$
\iota : P_{X_n,\lambda} \ni \beta^\vee \mapsto w_\lambda \beta^\vee \in P_{X_n,\lambda}
$$
gives an involutive anti-automorphism of the poset $P_{X_n,\lambda}$.
\end{prop}

\begin{demo}{Proof}
It is enough to show that $\beta^\vee > \gamma^\vee$ implies $w_\lambda \beta^\vee < w_\lambda \gamma^\vee$ 
for $\beta^\vee$, $\gamma^\vee \in P_{X_n,\lambda}$.
It follows from $\langle \beta^\vee, \lambda \rangle = \langle \gamma^\vee, \lambda \rangle = 1$ that 
$\beta^\vee - \gamma^\vee$ is a linear combination of $\Pi^\vee \setminus \{ \alpha^\vee \}$ 
with nonnegative integer coefficients, 
where $\Pi^\vee$ is the set of simple coroots and $\alpha^\vee$ is the simple coroot 
dual to $\lambda$.
Since $w_\lambda (\Pi^\vee \setminus \{ \alpha^\vee \}) = - (\Pi^\vee \setminus \{ \alpha^\vee \})$, 
we see that $w_\lambda \beta^\vee - w_\lambda \gamma^\vee$ 
is a linear combination of $\Pi^\vee \setminus \{ \alpha^\vee \}$ 
with nonpositive integer coefficients.
\qed
\end{demo}

The following properties of minuscule posets can be checked easily 
(e.g., by using a description given in the next subsection).

\begin{prop}
\label{prop:minuscule2}
Let $P = P_{X_n,\lambda}$ be the minuscule poset corresponding to a minuscule weight $\lambda$, 
and $c : P \to \Pi$ the coloring.
\begin{enumerate}
\item[(a)]
The poset $P$ is graded of height $h-1$, where $h$ is the Coxeter number of $\mathfrak{g}$.
\item[(b)]
The poset $P$ has a unique minimal element $v_{\min}$ and a unique maximal element $v_{\max}$.
Moreover, if we put $\alpha_{\min} = c(v_{\min})$ and $\alpha_{\max} = c(v_{\max})$, 
then the simple root $\alpha_{\min}$ corresponds to the fundamental weight $\lambda$ 
and $\alpha_{\max} = -w_0 \alpha_{\min}$ corresponds to $-w_0 \lambda$, 
where $w_0$ is the longest element of $W$.
\item[(c)]
If $v \lessdot w$ in $P$, then their colors $c(v)$ and $c(w)$ are adjacent in 
the Dynkin diagram of $\mathfrak{g}$.
\item[(d)]
For each $\alpha \in \Pi$, the subposet $P^\alpha = \{ v \in P : c(v) = \alpha \}$ is a chain.
\item[(e)]
If $v$, $w \in P^\alpha$, then the difference $\rank(v) - \rank(w)$ is even.
\end{enumerate}
\end{prop}

\subsection{%
Description of minuscule posets
}

In this subsection we give an explicit description of minuscule posets and their colorings.
The minuscule posets can be embedded into the poset $\Int^2$, 
where $(i,j) \le (i',j')$ in $\Int^2$ if and only if $i \le i'$ and $j \le j'$.

\begin{figure}[ht]
\centering
\begin{picture}(130,130)
\put(5,45){\circle{10}}
\put(25,25){\circle{10}}
\put(25,65){\circle{10}}
\put(45,5){\circle{10}}
\put(45,45){\circle{10}}
\put(45,85){\circle{10}}
\put(65,25){\circle{10}}
\put(65,65){\circle{10}}
\put(65,105){\circle{10}}
\put(85,45){\circle{10}}
\put(85,85){\circle{10}}
\put(85,125){\circle{10}}
\put(105,65){\circle{10}}
\put(105,105){\circle{10}}
\put(125,85){\circle{10}}
\put(25,25){\circle{10}}
\put(9,49){\line(1,1){12}}
\put(9,41){\line(1,-1){12}}
\put(29,29){\line(1,1){12}}
\put(29,21){\line(1,-1){12}}
\put(29,69){\line(1,1){12}}
\put(29,61){\line(1,-1){12}}
\put(49,9){\line(1,1){12}}
\put(49,49){\line(1,1){12}}
\put(49,41){\line(1,-1){12}}
\put(49,89){\line(1,1){12}}
\put(49,81){\line(1,-1){12}}
\put(69,29){\line(1,1){12}}
\put(69,69){\line(1,1){12}}
\put(69,61){\line(1,-1){12}}
\put(69,109){\line(1,1){12}}
\put(69,101){\line(1,-1){12}}
\put(89,49){\line(1,1){12}}
\put(89,89){\line(1,1){12}}
\put(89,81){\line(1,-1){12}}
\put(89,121){\line(1,-1){12}}
\put(109,69){\line(1,1){12}}
\put(109,101){\line(1,-1){12}}
\put(0,40){\makebox(10,10){$1$}}
\put(20,20){\makebox(10,10){$2$}}
\put(20,60){\makebox(10,10){$2$}}
\put(40,0){\makebox(10,10){$3$}}
\put(40,40){\makebox(10,10){$3$}}
\put(40,80){\makebox(10,10){$3$}}
\put(60,20){\makebox(10,10){$4$}}
\put(60,60){\makebox(10,10){$4$}}
\put(60,100){\makebox(10,10){$4$}}
\put(80,40){\makebox(10,10){$5$}}
\put(80,80){\makebox(10,10){$5$}}
\put(80,120){\makebox(10,10){$5$}}
\put(100,60){\makebox(10,10){$6$}}
\put(100,100){\makebox(10,10){$6$}}
\put(120,80){\makebox(10,10){$7$}}
\end{picture}
\caption{$P_{A_7,\varpi_3}$}
\label{fig:a}
\end{figure}
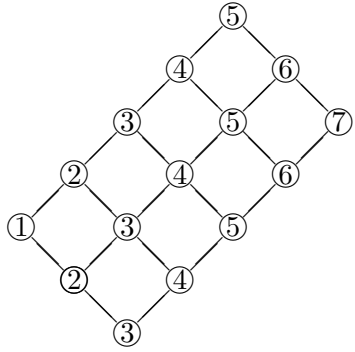

\paragraph{Type $A_n$.}
The positive coroot system $\Delta^\vee_+$ of type $A_n$ can be described as 
$\Delta^\vee_+ = \{ e_i - e_j : 1 \le i < j \le n+1 \}$ with $e_1 + \dots + e_{n+1} = 0$.
Then we have 
$$
P_{A_n,\varpi_r}
 =
\{ e_i - e_j : 1 \le i \le r, \ r+1 \le j \le n+1 \}
$$
and tha map $e_i - e_j \mapsto (r-i,j-r-1)$ gives an isomorphism of posets 
from $P_{A_n,\varpi_r}$ to the subposet 
$$
\{ (i,j) \in \Int^2 : 0 \le i \le r-1, \, 0 \le j \le n-r \}
 \subset \Int^2.
$$
The poset $P_{A_n,\varpi_r}$ is a product poset $[0,r-1] \times [0,n-r]$ of two chains, 
where $[0,m] = \{ 0, 1, \dots, m \}$ is a chain.
We call this poset $P_{A_n,\varpi_r}$ a \emph{rectangle poset}.
The involution $\iota$ is the $180^\circ$ rotation of the Hasse diagram.
For example, the Hasse diagram and the coloring of $P_{A_7, \varpi_3}$ are given 
in Figure~\ref{fig:a}, 
where we label a vertex $v$ with $i$ to indicate that $c(v) = \alpha_i$. 

\paragraph{Type $B_n$.}
If we realize the positive coroot system $\Delta^\vee_+$ of type $B_n$ as 
$\Delta^\vee_+ = \{ e_i \pm e_j : 1 \le i < j \le n \} \cup \{ 2 e_i : 1 \le i \le n \}$, then 
we have
$$
P_{B_n,\varpi_n}
 =
\{ e_i + e_j : 1 \le i \le j \le n \},
$$
and the map $e_i + e_j \mapsto (n-j,n-i)$ gives an poset isomorphism from $P_{B_n,\varpi_n}$ 
to the subposet 
$$
\{ (i,j) \in \Int^2 : 0 \le i \le j \le n-1 \}
 \subset \Int^2.
$$
We call $P_{B_n,\varpi_n}$ a \emph{shifted staircase poset}.
The involution $\iota$ is the horizontal flip of the Hasse diagram.
For example the Hasse diagram of $P_{B_4,\varpi_4}$ and its coloring are given in Figure~\ref{fig:b}.
\begin{figure}[ht]
\setlength{\unitlength}{1pt}
\begin{minipage}{0.24\hsize}
\centering
\begin{picture}(70,130)
\put(5,5){\circle{10}}
\put(5,45){\circle{10}}
\put(5,85){\circle{10}}
\put(5,125){\circle{10}}
\put(25,25){\circle{10}}
\put(25,65){\circle{10}}
\put(25,105){\circle{10}}
\put(45,45){\circle{10}}
\put(45,85){\circle{10}}
\put(65,65){\circle{10}}
\put(9,9){\line(1,1){12}}
\put(9,49){\line(1,1){12}}
\put(9,41){\line(1,-1){12}}
\put(9,89){\line(1,1){12}}
\put(9,81){\line(1,-1){12}}
\put(9,121){\line(1,-1){12}}
\put(29,29){\line(1,1){12}}
\put(29,69){\line(1,1){12}}
\put(29,61){\line(1,-1){12}}
\put(29,101){\line(1,-1){12}}
\put(49,49){\line(1,1){12}}
\put(49,81){\line(1,-1){12}}
\put(0,0){\makebox(10,10){$4$}}
\put(0,40){\makebox(10,10){$4$}}
\put(0,80){\makebox(10,10){$4$}}
\put(0,120){\makebox(10,10){$4$}}
\put(20,20){\makebox(10,10){$3$}}
\put(20,60){\makebox(10,10){$3$}}
\put(20,100){\makebox(10,10){$3$}}
\put(40,40){\makebox(10,10){$2$}}
\put(40,80){\makebox(10,10){$2$}}
\put(60,60){\makebox(10,10){$1$}}
\end{picture}
\caption{$P_{B_4,\varpi_4}$}
\label{fig:b}
\end{minipage}
\begin{minipage}{0.24\hsize}
\centering
\begin{picture}(70,130)
\put(5,5){\circle{10}}
\put(5,125){\circle{10}}
\put(25,25){\circle{10}}
\put(25,105){\circle{10}}
\put(45,45){\circle{10}}
\put(45,85){\circle{10}}
\put(65,65){\circle{10}}
\put(9,9){\line(1,1){12}}
\put(9,121){\line(1,-1){12}}
\put(29,29){\line(1,1){12}}
\put(29,101){\line(1,-1){12}}
\put(49,49){\line(1,1){12}}
\put(49,81){\line(1,-1){12}}
\put(0,0){\makebox(10,10){$1$}}
\put(0,120){\makebox(10,10){$1$}}
\put(20,20){\makebox(10,10){$2$}}
\put(20,100){\makebox(10,10){$2$}}
\put(40,40){\makebox(10,10){$3$}}
\put(40,80){\makebox(10,10){$3$}}
\put(60,60){\makebox(10,10){$4$}}
\end{picture}
\caption{$P_{C_4,\varpi_1}$}
\label{fig:c}
\end{minipage}
\begin{minipage}{0.24\hsize}
\centering
\begin{picture}(70,130)
\put(5,5){\circle{10}}
\put(5,125){\circle{10}}
\put(25,25){\circle{10}}
\put(25,65){\circle{10}}
\put(25,105){\circle{10}}
\put(45,45){\circle{10}}
\put(45,85){\circle{10}}
\put(65,65){\circle{10}}
\put(9,9){\line(1,1){12}}
\put(9,121){\line(1,-1){12}}
\put(29,29){\line(1,1){12}}
\put(29,69){\line(1,1){12}}
\put(29,61){\line(1,-1){12}}
\put(29,101){\line(1,-1){12}}
\put(49,49){\line(1,1){12}}
\put(49,81){\line(1,-1){12}}
\put(0,0){\makebox(10,10){$1$}}
\put(0,120){\makebox(10,10){$1$}}
\put(20,20){\makebox(10,10){$2$}}
\put(20,60){\makebox(10,10){$5$}}
\put(20,100){\makebox(10,10){$2$}}
\put(40,40){\makebox(10,10){$3$}}
\put(40,80){\makebox(10,10){$3$}}
\put(60,60){\makebox(10,10){$4$}}
\end{picture}
\caption{$P_{D_5,\varpi_1}$}
\label{fig:d1}
\end{minipage}
\begin{minipage}{0.24\hsize}
\centering
\begin{picture}(70,130)
\put(5,5){\circle{10}}
\put(5,45){\circle{10}}
\put(5,85){\circle{10}}
\put(5,125){\circle{10}}
\put(25,25){\circle{10}}
\put(25,65){\circle{10}}
\put(25,105){\circle{10}}
\put(45,45){\circle{10}}
\put(45,85){\circle{10}}
\put(65,65){\circle{10}}
\put(9,9){\line(1,1){12}}
\put(9,49){\line(1,1){12}}
\put(9,41){\line(1,-1){12}}
\put(9,89){\line(1,1){12}}
\put(9,81){\line(1,-1){12}}
\put(9,121){\line(1,-1){12}}
\put(29,29){\line(1,1){12}}
\put(29,69){\line(1,1){12}}
\put(29,61){\line(1,-1){12}}
\put(29,101){\line(1,-1){12}}
\put(49,49){\line(1,1){12}}
\put(49,81){\line(1,-1){12}}
\put(0,0){\makebox(10,10){$5$}}
\put(0,40){\makebox(10,10){$4$}}
\put(0,80){\makebox(10,10){$5$}}
\put(0,120){\makebox(10,10){$4$}}
\put(20,20){\makebox(10,10){$3$}}
\put(20,60){\makebox(10,10){$3$}}
\put(20,100){\makebox(10,10){$3$}}
\put(40,40){\makebox(10,10){$2$}}
\put(40,80){\makebox(10,10){$2$}}
\put(60,60){\makebox(10,10){$1$}}
\end{picture}
\caption{$P_{D_5,\varpi_5}$}
\label{fig:d2}
\end{minipage}
\end{figure}

\paragraph{Type $C_n$.}
If we realize the positive coroot system $\Delta^\vee_+$ of type $C_n$ as 
$\Delta^\vee_+ = \{ e_i \pm e_j : 1 \le i < j \le n \} \cup \{ e_i : 1 \le i \le n \}$, then 
we have
$$
P_{C_n,\varpi_1}
 =
\{ e_1 - e_2, \dots, e_1 - e_n, e_1, e_1 + e_n, \dots, e_1 + e_2 \}.
$$
The poset $P_{C_n,\varpi_1}$ is a chain, and isomorphic to 
the subposet 
$$
\{ (1,1), \dots, (1,n-1), (1,n), (2,n), \dots, (n,n) \} \subset \Int^2.
$$
For example the Hasse diagram of $P_{C_4,\varpi_1}$ and its coloring are given in Figure~\ref{fig:c}.
Note that $P_{C_n,\varpi_1}$ is isomorphic to $P_{A_{2n-1},\varpi_1}$,
but they have different colorings.

\paragraph{Type $D_n$.}
We realize the positive coroot system $\Delta^\vee_+$ of type $D_n$ as 
$\Delta^\vee_+ = \{ e_i \pm e_j : 1 \le i < j \le n \}$.

For the minuscule weight $\varpi_1$, we have
$$
P_{D_n,\varpi_1}
 =
\{ e_1-e_2, \dots, e_1-e_{n-1}, e_1-e_n, e_1+e_n, e_1+e_{n-1}, \dots, e_1+e_2 \},
$$
and it is isomorphic to the subposet
$$
\{ (1,1), \dots, (1,n-1), (1,n), (2,n-1), (2,n), \dots, (n,n) \} \subset \Int^2.
$$
See Figure~\ref{fig:d1} for the Hasse diagram of $P_{D_5,\varpi_1}$ and its coloring. 
The poset $P_{D_n,\varpi_1}$ is called a \emph{double-tailed diamond poset}. 
The involutive anti-automorphism $\iota$ is given by
$$
\iota (e_1 + e_k) = e_1 - e_k \quad(1 \le k \le n-1),
\quad
\iota (e_1 + \ep e_n) = e_1 + (-1)^n \ep e_n.
$$

For the minuscule weights $\varpi_n$ and $\varpi_{n-1}$, we have 
$$
P_{D_n,\varpi_n}
 =
\{ e_i + e_j : 1 \le i < j \le n \}
$$
and $P_{D_n,\varpi_{n-1}}$ is obtained from $P_{D_n,\varpi_n}$ by replacing $e_i + e_n$ with $e_i-e_n$ 
for $1 \le i \le n-1$.
Both posets $P_{D_n,\varpi_n}$ and $P_{D_n,\varpi_{n-1}}$ are isomorphic to 
$\{ (i,j) \in \Int^2 : 0 \le i \le j \le n-2 \}$.
For example, the Hasse diagram and the coloring of $P_{D_5,\varpi_5}$ are given in Figure~\ref{fig:d2}.
Note that $P_{D_n,\varpi_{n-1}} \cong P_{D_n,\varpi_n}$ and they are isomorphic to $P_{B_{n-1},\varpi_{n-1}}$, 
but they have different colorings.

\paragraph{Type $E_6$.}
The minuscule poset $P_{E_6,\varpi_6}$ is isomorphic to the subposet
$$
\left\{
\begin{matrix}
 (1,1),(2,1),(3,1),(4,1),(5,1),
             (3,2),(4,2),(5,2), \\
                   (4,3),(5,3),(6,3),
                   (4,4),(5,4),(6,4),(7,4),(8,4)
\end{matrix}
\right\}
\subset \Int^2,
$$
and the Hasse diagram and the coloring are given in Figure~\ref{fig:e6}.
The involution $\iota$ is the $180^\circ$ rotation of the Hasse diagram.
As posets, $P_{E_6,\varpi_1} \cong P_{E_6,\varpi_6}$.
\begin{figure}[ht]
\begin{minipage}{0.49\hsize}
\begin{center}
\begin{picture}(90,210)
\put(5,85){\circle{10}}
\put(5,205){\circle{10}}
\put(25,65){\circle{10}}
\put(25,105){\circle{10}}
\put(25,145){\circle{10}}
\put(25,185){\circle{10}}
\put(45,45){\circle{10}}
\put(45,85){\circle{10}}
\put(45,125){\circle{10}}
\put(45,165){\circle{10}}
\put(65,25){\circle{10}}
\put(65,65){\circle{10}}
\put(65,105){\circle{10}}
\put(65,145){\circle{10}}
\put(85,5){\circle{10}}
\put(85,125){\circle{10}}
\put(9,89){\line(1,1){12}}
\put(9,81){\line(1,-1){12}}
\put(9,201){\line(1,-1){12}}
\put(29,69){\line(1,1){12}}
\put(29,61){\line(1,-1){12}}
\put(29,109){\line(1,1){12}}
\put(29,101){\line(1,-1){12}}
\put(29,149){\line(1,1){12}}
\put(29,141){\line(1,-1){12}}
\put(29,181){\line(1,-1){12}}
\put(49,49){\line(1,1){12}}
\put(49,41){\line(1,-1){12}}
\put(49,89){\line(1,1){12}}
\put(49,81){\line(1,-1){12}}
\put(49,129){\line(1,1){12}}
\put(49,121){\line(1,-1){12}}
\put(49,161){\line(1,-1){12}}
\put(69,21){\line(1,-1){12}}
\put(69,109){\line(1,1){12}}
\put(69,141){\line(1,-1){12}}
\put(0,80){\makebox(10,10){$1$}}
\put(0,200){\makebox(10,10){$1$}}
\put(20,60){\makebox(10,10){$3$}}
\put(20,100){\makebox(10,10){$3$}}
\put(20,140){\makebox(10,10){$2$}}
\put(20,180){\makebox(10,10){$3$}}
\put(40,40){\makebox(10,10){$4$}}
\put(40,80){\makebox(10,10){$4$}}
\put(40,120){\makebox(10,10){$4$}}
\put(40,160){\makebox(10,10){$4$}}
\put(60,20){\makebox(10,10){$5$}}
\put(60,60){\makebox(10,10){$2$}}
\put(60,100){\makebox(10,10){$5$}}
\put(60,140){\makebox(10,10){$5$}}
\put(80,0){\makebox(10,10){$6$}}
\put(80,120){\makebox(10,10){$6$}}
\end{picture}
\end{center}
\caption{$P_{E_6,\varpi_6}$}
\label{fig:e6}
\end{minipage}
\begin{minipage}{0.49\hsize}
\begin{center}
\begin{picture}(110,330)
\put(5,5){\circle{10}}
\put(5,165){\circle{10}}
\put(5,325){\circle{10}}
\put(25,25){\circle{10}}
\put(25,145){\circle{10}}
\put(25,185){\circle{10}}
\put(25,305){\circle{10}}
\put(45,45){\circle{10}}
\put(45,85){\circle{10}}
\put(45,125){\circle{10}}
\put(45,165){\circle{10}}
\put(45,205){\circle{10}}
\put(45,245){\circle{10}}
\put(45,285){\circle{10}}
\put(65,65){\circle{10}}
\put(65,105){\circle{10}}
\put(65,145){\circle{10}}
\put(65,185){\circle{10}}
\put(65,225){\circle{10}}
\put(65,265){\circle{10}}
\put(85,85){\circle{10}}
\put(85,125){\circle{10}}
\put(85,165){\circle{10}}
\put(85,205){\circle{10}}
\put(85,245){\circle{10}}
\put(105,105){\circle{10}}
\put(105,225){\circle{10}}
\put(9,9){\line(1,1){12}}
\put(9,169){\line(1,1){12}}
\put(9,161){\line(1,-1){12}}
\put(9,321){\line(1,-1){12}}
\put(29,29){\line(1,1){12}}
\put(29,149){\line(1,1){12}}
\put(29,141){\line(1,-1){12}}
\put(29,189){\line(1,1){12}}
\put(29,181){\line(1,-1){12}}
\put(29,301){\line(1,-1){12}}
\put(49,49){\line(1,1){12}}
\put(49,89){\line(1,1){12}}
\put(49,81){\line(1,-1){12}}
\put(49,129){\line(1,1){12}}
\put(49,121){\line(1,-1){12}}
\put(49,169){\line(1,1){12}}
\put(49,161){\line(1,-1){12}}
\put(49,209){\line(1,1){12}}
\put(49,201){\line(1,-1){12}}
\put(49,249){\line(1,1){12}}
\put(49,241){\line(1,-1){12}}
\put(49,281){\line(1,-1){12}}
\put(69,69){\line(1,1){12}}
\put(69,109){\line(1,1){12}}
\put(69,101){\line(1,-1){12}}
\put(69,149){\line(1,1){12}}
\put(69,141){\line(1,-1){12}}
\put(69,189){\line(1,1){12}}
\put(69,181){\line(1,-1){12}}
\put(69,229){\line(1,1){12}}
\put(69,221){\line(1,-1){12}}
\put(69,261){\line(1,-1){12}}
\put(89,89){\line(1,1){12}}
\put(89,121){\line(1,-1){12}}
\put(89,209){\line(1,1){12}}
\put(89,241){\line(1,-1){12}}
\put(0,0){\makebox(10,10){$7$}}
\put(0,160){\makebox(10,10){$7$}}
\put(0,320){\makebox(10,10){$7$}}
\put(20,20){\makebox(10,10){$6$}}
\put(20,140){\makebox(10,10){$6$}}
\put(20,180){\makebox(10,10){$6$}}
\put(20,300){\makebox(10,10){$6$}}
\put(40,40){\makebox(10,10){$5$}}
\put(40,80){\makebox(10,10){$2$}}
\put(40,120){\makebox(10,10){$5$}}
\put(40,160){\makebox(10,10){$5$}}
\put(40,200){\makebox(10,10){$5$}}
\put(40,240){\makebox(10,10){$2$}}
\put(40,280){\makebox(10,10){$5$}}
\put(60,60){\makebox(10,10){$4$}}
\put(60,100){\makebox(10,10){$4$}}
\put(60,140){\makebox(10,10){$4$}}
\put(60,180){\makebox(10,10){$4$}}
\put(60,220){\makebox(10,10){$4$}}
\put(60,260){\makebox(10,10){$4$}}
\put(80,80){\makebox(10,10){$3$}}
\put(80,120){\makebox(10,10){$3$}}
\put(80,160){\makebox(10,10){$2$}}
\put(80,200){\makebox(10,10){$3$}}
\put(80,240){\makebox(10,10){$3$}}
\put(100,100){\makebox(10,10){$1$}}
\put(100,220){\makebox(10,10){$1$}}
\end{picture}
\end{center}
\caption{$P_{E_7,\varpi_7}$}
\label{fig:e7}
\end{minipage}
\end{figure}

\paragraph{Type $E_7$.}
The minuscule poset $P_{E_7,\varpi_7}$ is isomorphic to the subposet
$$
\left\{
\begin{matrix}
     (1,1),(1,2),(1,3),(1,4),(1,5),(1,6),
                       (2,4),(2,5),(2,6), \\
                             (3,5),(3,6),(3,7),
                             (4,5),(4,6),(4,7),
                             (5,5),(5,6),(5,7), \\
(4,8),(4,9),
(5,8),(5,9),
                                               (6,8),(6,9),
                                                     (7,9),
                                                     (8,9),
                                                     (9,9)
\end{matrix}
\right\}
\subset \Int^2,
$$
and the Hasse diagram and the coloring are given in Figure~\ref{fig:e7}.
The involution $\iota$ is the horizontal flip of the Hasse diagram.

\section{%
Periodicity
}

The goal of this section is to prove the periodicity of birational rowmotion 
and Coxeter-motion (Theorem~\ref{thm:main1} (a) and Theorem~\ref{thm:main2} (a)).

\subsection{%
Periodicity of birational rowmotion
}

For the birational rowmotion map on minuscule posets, 
the periodicity has been established in \cite{GR1, GR2} 
except for the type $E_7$ minuscule poset.
Let $P$ be a minuscule poset associated to a Lie algebra $\mathfrak{g}$, 
and $\rho^{A,B} : \KK^{A,B}(P) \to \KK^{A,B}$ the birational rowmotion map.
Since the periodicity depends only on the poset structure, 
we may assume that $\mathfrak{g}$ is simply-laced.
And by Proposition~\ref{prop:minuscule2} (a), Lemmas~\ref{lem:order} and \ref{lem:A=B=1}, 
it is enough to show that 
$\rho = \rho^{1,1}$ satisfies $\rho^h = 1$, 
where $h$ is the Coxeter number of $\mathfrak{g}$.
\begin{itemize}
\item
If $P$ is a type $A_n$ minuscule poset, i.e., if $P$ is a rectangle poset $[0,r-1] \times [0, n-r]$,
then it was shown that the birational rowmotion map $\rho$ has order $n+1$ 
(Grinberg--Roby \cite[Theorem~30]{GR2}, see \cite[Corollary~2.12]{MR} for another proof).
\item
If $P = P_{D_n,\varpi_1}$ is a double-tailed diamond poset, 
then $P$ is a skeletal poset of height $2n-3$, 
and it follows from \cite[Propositions~61, 74 and 75]{GR1} that $\rho$ has order $2n-2$ 
(see \cite[Section~10]{GR1} for a definition of skeletal posets and details).
\item
If $P = P_{D_n,\varpi_n}$ is a shifted staircase poset, 
then Grinberg--Roby \cite[Theorem58]{GR2} proved that $\rho$ has order $2n$.
\item
If $P = P_{E_6,\varpi_6}$ is the minuscule poset of type $E_6$, then by using a computer we can verify 
that $\rho$ has order $12$. 
\item
Let $P = P_{E_7,\varpi_7}$ be the minuscule poset of type $E_7$.
Given an initial state $X \in \KK^{1,1}(P)$, 
we regard $\{ X(v) : v \in P \}$ as indeterminates 
and introduce new indeterminates $\{ Z(v) : v \in P \}$ by (\ref{eq:X2Z}).
With the author's laptop, it takes about 20 seconds for Maple19 
to compute all the values $\left( \rho^k X \right)(v)$ ($0 \le k \le 18$, $v \in P$) 
as rational functions in $\{ Z(v) : v \in P \}$ 
and check that $\left( \rho^{18} X \right)(v) = X(v)$ for all $v \in P$.
\end{itemize}
This completes the proof of Theorem~\ref{thm:main1} (a).

\subsection{%
Periodicity of birational Coxeter-motion
}

In order to prove the periodicity of birational Coxeter-motion (Theorem~\ref{thm:main2} (a)), 
we work with the birational toggle group and show that 
any birational Coxeter-motion maps are conjugate to the birational rowmotion map 
in this group.

Let $P$ be a finite poset and fix positive real numbers $A$ and $B$.
We define the \emph{birational toggle group}, denote by $G(P)$, 
to be the subgroup generated by birational toggles $\tau_v = \tau^{A,B}_v$ ($v \in P$) 
in the group of all bijections on $\KK^{A,B}(P)$.

A key tool here is the non-commutativity graph.
Given elements $g_1, \dots, g_n$ of a group $G$, 
the \emph{non-commutativity graph} $\Gamma(g_1, \dots, g_n)$ is defined as the graph 
with vertex set $\{ 1, 2, \dots, n \}$, in which 
two vertices $i$ and $j$ are joined if and only if $g_i g_j \neq g_j g_i$.
The following lemma is useful.

\begin{lemma}
\label{lem:conj}
(\cite[V, \S6, n${}^\circ$1, Lemma~1]{B1})
Let $g_1, \dots, g_n$ be elements of a group $G$.
If the non-commutativity graph $\Gamma(g_1, \dots, g_n)$ has no cycle, 
then $g_{\nu(1)} \dots g_{\nu(n)}$ is conjugate to $g_1 \dots g_n$ in $G$
for any permutation $\nu \in S_n$.
\end{lemma}

First we prove that all birational Coxeter-motion maps are conjugate.

\begin{prop}
\label{prop:conj}
Let $P$ be a minuscule poset.
Then all birational Coxeter-motion maps are conjugate to each other in 
the birational toggle group $G(P)$.
\end{prop}

\begin{demo}{Proof}
Note that birational toggles $\tau_v$ and $\tau_w$ are commutative 
unless $v \lessdot w$ ore $v \gtrdot w$.
It follows from Proposition~\ref{prop:minuscule2} (c) that,  
if simple roots $\alpha$ and $\beta$ are not adjacent in the Dynkin diagram of $\mathfrak{g}$, 
then the corresponding elements $\sigma_\alpha$ and $\sigma_\beta$ commute with each other in $G(P)$.
Hence the non-commutativity graph $\Gamma(\sigma_{\alpha_1}, \dots, \sigma_{\alpha_n})$, 
where $\alpha_1, \dots, \alpha_n$ are the simple roots, 
is a subgraph (of the underlying simple graph) of the Dynkin diagram.
Since the Dynkin diagram of $\mathfrak{g}$ has no cycle, 
we can use Lemma~\ref{lem:conj} to conclude that any two Coxeter-motion maps 
are conjugate in $G(P)$.
\qed
\end{demo}

The periodicity of birational Coxeter-motion maps (Theorem~\ref{thm:main2} (a)) immediately follows 
from the following thoerem and the periodicity of the birational rowmotion map 
(Theorem~\ref{thm:main1} (a)).

\begin{theorem}
\label{thm:conj}
Let $P$ be a minuscule poset.
Then any birational Coxeter-motion map is conjugate to the birational rowmotion map $\rho = \rho^{A,B}$ 
in the biratoinal toggle group $G(P)$.
\end{theorem}

This theorem is a birational lift of \cite[Theorem~1.3]{RS}.
In order to prove this theorem, we use a notion of rc-poset, 
which is introduced by Striker--Williams \cite[Section~4.2]{SW}.
We put $\Lambda = \{ (i,j) \in \Int^2 : \text{$i+j$ is even} \}$.
A poset $P$ is called a \emph{rowed-and-columned poset} (\emph{rc-poset} for short) 
if there is a map $\pi : P \to \Lambda$ such that, 
if $v$ covers $u$ in $P$ and $\pi(v) = (i,j)$, 
then $\pi(u) = (i+1,j-1)$ or $(i-1,j-1)$.
Minuscule posets $P = P_{X_n,\lambda}$ are rc-posets with respect to the composition map 
$\pi : P \to \Lambda$ of the embedding $P \hookrightarrow \Int^2$ given in Subsection~3.2 
and the map $\Int^2 \ni (i,j) \mapsto (j-i,j+i) \in \Lambda$.
A \emph{row} (resp. \emph{column}) of an RC-poset $P$ is a subset $M$ of $P$ of the form 
\begin{gather*}
M = \{ v \in P : \text{the second coordinate of $\pi(v)$ equals $r$} \},
\\
\text{(resp. }
M = \{ v \in P : \text{the first coordinate of $\pi(v)$ equals $c$} \}
\text{)}
\end{gather*}
for some $r$ (resp. $c$).
If $M$ is a subset of a row or a column of $P$, 
then the composition of toggles $\tau_v$ ($v \in M$) is independent of the order of composition, 
so we denote by $\tau[M]$ the resulting element of the toggle group $G(P)$.
If $R_1, \dots, R_n$ are the non-empty rows of an rc-poset $P$ from bottom to top, 
then the rowmotion map $\rho = \rho^{A,B}$ is given by
$$
\rho = \tau[R_1] \circ \tau[R_2] \circ \dots \circ \tau[R_n].
$$
The following Lemma is proved by exactly the same argument as in \cite{SW}.

\begin{lemma}
\label{lem:rc}
(\cite[Theorem~5.2]{SW})
Let $P$ be an rc-poset.
Let $R_1, \dots, R_n$ be the non-empty rows of $P$ from bottom to top, and 
$C_1, \dots, C_m$ the non-empty columns of $P$ from left to right.
Then the rowmotion map $\rho$ is conjugate to 
$\tau[C_{\nu(1)}] \circ \cdots \circ \tau[C_{\nu(m)}]$ in $G(P)$ for any $\nu \in S_m$.
\end{lemma}

We prove Theorem~\ref{thm:conj} by using this lemma.

\begin{demo}{Proof of Theorem~\ref{thm:conj}}
Let $\Pi = \{ \alpha_1, \dots, \alpha_n \}$ be the set of simple roots, 
where we follow the numbering in \cite{B1}, 
and $C_1, \dots, C_m$ the non-empty columns of $P$ 
(see Figures~\ref{fig:a}--\ref{fig:e7}).
Then, by Lemmas~\ref{lem:conj} and \ref{lem:rc}, 
it is enough to prove that $\gamma = \sigma_{\alpha_1} \cdots \sigma_{\alpha_n}$ is conjugate to 
$\tau[C_1] \cdots \tau[C_m]$.
We prove this claim by a case-by-case argument.
\begin{itemize}
\item
If $P = P_{A_n,\varpi_r}$, 
then $\sigma_{\alpha_i} = \tau[C_i]$ for $1 \le i \le n$ 
and $\gamma = \tau[C_1] \cdots \tau[C_m]$.
\item
If $P = P_{B_n,\varpi_n}$, 
then $\sigma_{\alpha_i} = \tau[C_{n+1-i}]$ for $1 \le i \le n$, 
and $\gamma$ is conjugate to $\sigma_{\alpha_n} \cdots \sigma_{\alpha_1}
 = \tau[C_1] \cdots \tau[C_m]$ by Lemma~\ref{lem:conj}.
\item
If $P = P_{C_n,\varpi_1}$, 
then $\sigma_{\alpha_i} = \tau[C_i]$ for $1 \le i \le n$ 
and $\gamma = \tau[C_1] \cdots \tau[C_n]$.
\item
If $P = P_{D_n,\varpi_1}$, 
then $\tau[C_i] = \sigma_{\alpha_i}$ for $i \neq n-3$ 
and $\tau[C_{n-3}] = \sigma_{\alpha_{n-3}} \sigma_{\alpha_n}$.
Hence $\tau[C_1] \cdots \tau[C_{n-1}] = \sigma_{\alpha_1} \cdots \sigma_{\alpha_{n-4}}
 \sigma_{\alpha_{n-3}} \sigma_{\alpha_n} \sigma_{\alpha_{n-2}} \sigma_{\alpha_{n-1}}$ 
is conjugate to $\gamma$ by Lemma~\ref{lem:conj}.
\item
If $P = P_{D_n, \varpi_n}$, 
then $\tau[C_1] = \sigma_{\alpha_{n-1}} \sigma_{\alpha_n}$ 
and $\tau[C_i] = \sigma_{\alpha_{n-i}}$ for $2 \le i \le n-1$.
Hence $\tau[C_1] \cdots \tau[C_{n-1}] = \sigma_{\alpha_{n-1}} \sigma_{\alpha_n} 
\sigma_{\alpha_{n-2}} \cdots \sigma_{\alpha_1}$ 
is conjugate to $\gamma$ by Lemma~\ref{lem:conj}.
\item
If $P = P_{E_6,\varpi_6}$, then we have
$$
C_1 = P^{\alpha_1},
\quad
C_2 = P^{\alpha_3} \sqcup (C_2 \cap P^{\alpha_2}),
\quad
C_3 = P^{\alpha_4},
\quad
C_4 = P^{\alpha_5} \sqcup (C_4 \cap P^{\alpha_2}),
\quad
C_5 = P^{\alpha_6}.
$$
If we put
\begin{gather*}
g_1 = \tau[P^{\alpha_1}],
\quad
g_2 = \tau[P^{\alpha_3}]
\quad
g_3 = \tau[P^{\alpha_4}],
\quad
g_4 = \tau[P^{\alpha_5}],
\quad
g_5 = \tau[P^{\alpha_6}],
\\
g_6 = \tau[C_2 \cap P^{\alpha_2}],
\quad
g_7 = \tau[C_5 \cap P^{\alpha_2}],
\end{gather*}
then Figure~\ref{fig:e6-nc} shows the non-commutativity graph $\Gamma(g_1, \dots, g_7)$.
\begin{figure}
\centering
\setlength{\unitlength}{1pt}
\begin{picture}(130,30)
\put(5,25){\circle{10}}
\put(35,25){\circle{10}}
\put(65,25){\circle{10}}
\put(95,25){\circle{10}}
\put(125,25){\circle{10}}
\put(45,5){\circle{10}}
\put(85,5){\circle{10}}
\put(10,25){\line(1,0){20}}
\put(40,25){\line(1,0){20}}
\put(70,25){\line(1,0){20}}
\put(100,25){\line(1,0){20}}
\put(49,9){\line(1,1){12}}
\put(81,9){\line(-1,1){12}}
\put(0,20){\makebox(10,10){$1$}}
\put(30,20){\makebox(10,10){$2$}}
\put(60,20){\makebox(10,10){$3$}}
\put(90,20){\makebox(10,10){$4$}}
\put(120,20){\makebox(10,10){$5$}}
\put(40,0){\makebox(10,10){$6$}}
\put(80,0){\makebox(10,10){$7$}}
\end{picture}
\caption{Non-commutativity graph for $P_{E_6,\varpi_6}$}
\label{fig:e6-nc}
\end{figure}
Hence by applying Lemma~\ref{lem:conj}, we see that 
$$
\tau[C_1] \cdots \tau[C_5]
 =
\tau[P^{\alpha_1}]
\tau[P^{\alpha_3}]
\tau[C_2 \cap P^{\alpha_2}]
\tau[P^{\alpha_4}]
\tau[P^{\alpha_5}]
\tau[C_5 \cap P^{\alpha_2}]
\tau[P^{\alpha_6}]
$$
is conjugate to
$$
\gamma =
\tau[P^{\alpha_1}]
\tau[C_2 \cap P^{\alpha_2}]
\tau[C_5 \cap P^{\alpha_2}]
\tau[P^{\alpha_3}]
\tau[P^{\alpha_4}]
\tau[P^{\alpha_5}]
\tau[P^{\alpha_6}].
$$
\item
If $P = P_{E_7,\varpi_7}$, then we have
\begin{gather*}
C_1 = P^{\alpha_7},
\quad
C_2 = P^{\alpha_6},
\quad
C_3 = (C_3 \cap P^{\alpha_2}) \sqcup P^{\alpha_5},
\\
C_4 = P^{\alpha_4},
\quad
C_5 = (C_5 \cap P^{\alpha_2}) \sqcup P^{\alpha_3},
\quad
C_6 = P^{\alpha_1},
\end{gather*}
and $P^{\alpha_2} = (C_3 \cap P^{\alpha_2}) \sqcup (C_5 \cap P^{\alpha_2})$.
If we put
\begin{gather*}
g_1 = \tau[P^{\alpha_7}],
\quad
g_2 = \tau[P^{\alpha_6}],
\quad
g_3 = \tau[P^{\alpha_5}],
\quad
g_4 = \tau[P^{\alpha_4}],
\quad
g_5 = \tau[P^{\alpha_3}],
\quad
g_6 = \tau[P^{\alpha_1}],
\\
g_7 = \tau[C_3 \cap P^{\alpha_2}],
\quad
g_8 = \tau[C_5 \cap P^{\alpha_2}]
\end{gather*}
then Figure~\ref{fig:e7-nc} shows the non-commutativity graph $\Gamma(g_1, \dots, g_8)$.
\begin{figure}
\centering
\setlength{\unitlength}{1pt}
\begin{picture}(160,30)
\put(5,25){\circle{10}}
\put(35,25){\circle{10}}
\put(65,25){\circle{10}}
\put(95,25){\circle{10}}
\put(125,25){\circle{10}}
\put(155,25){\circle{10}}
\put(75,5){\circle{10}}
\put(115,5){\circle{10}}
\put(10,25){\line(1,0){20}}
\put(40,25){\line(1,0){20}}
\put(70,25){\line(1,0){20}}
\put(100,25){\line(1,0){20}}
\put(130,25){\line(1,0){20}}
\put(79,9){\line(1,1){12}}
\put(111,9){\line(-1,1){12}}
\put(0,20){\makebox(10,10){$1$}}
\put(30,20){\makebox(10,10){$2$}}
\put(60,20){\makebox(10,10){$3$}}
\put(90,20){\makebox(10,10){$4$}}
\put(120,20){\makebox(10,10){$5$}}
\put(150,20){\makebox(10,10){$6$}}
\put(70,0){\makebox(10,10){$7$}}
\put(110,0){\makebox(10,10){$8$}}
\end{picture}
\caption{Non-commutativity graph for $P_{E_7,\varpi_7}$}
\label{fig:e7-nc}
\end{figure}
Hence by applying Lemma~\ref{lem:conj}, we see that 
$$
\tau[C_1] \cdots \tau[C_6]
 =
\tau[P^{\alpha_7}]
\tau[P^{\alpha_6}]
\tau[C_3 \cap P^{\alpha_2}] \tau[P^{\alpha_5}]
\tau[P^{\alpha_4}]
\tau[C_5 \cap P^{\alpha_2}] \tau[P^{\alpha_3}]
\tau[P^{\alpha_1}]
$$
is conjugate to 
$$
\gamma
 =
\tau[P^{\alpha_1}]
\tau[C_3 \cap P^{\alpha_2}]
\tau[C_5 \cap P^{\alpha_2}]
\tau[P^{\alpha_3}]
\tau[P^{\alpha_4}]
\tau[P^{\alpha_5}]
\tau[P^{\alpha_6}]
\tau[P^{\alpha_7}].
$$
\end{itemize}
This completes the proof of Theorem~\ref{thm:conj}, and hence of Theorem~\ref{thm:main2} (a).
\qed
\end{demo}

\section{%
Reciprocity
}

In this section we prove the reciprocity for birational rowmotion (Theorem~\ref{thm:main1} (b)) 
and propose a conjectural reciprocity for a particular birational Coxeter-motion map.

The proof of the reciprocity for birational rowmotion is based on a case-by-case analysis.
Let $P$ be a minuscule poset associated to a simple Lie algebra $\mathfrak{g}$ 
and $\rho^{A,B}$ the birational rowmotion map.
We may assume that $\mathfrak{g}$ is simply-laced and that $A=B=1$ (see Lemma~\ref{lem:A=B=1}).
For a type $A$ minuscule poset, 
the reciprocity was proved by Grinberg--Roby \cite[Theorem~30]{GR2} 
and Musiker--Roby \cite[Corollary~2.13]{MR}.
Also we can verify the reciprocity for the minuscule posets of types $E_6$ and $E_7$ 
by using a computer.
The remaining minuscule posets are the shifted staircase posets $P_{D_n,\varpi_n}$ 
and the double-tailed diamond posets $P_{D_n,\varpi_1}$.

\subsection{
Shifted staircase posets
}

Let $P = \{ (i,j) \in \Int^2 : 0 \le i \le j \le r \}$ be a shifted staircase poset, 
and $\rho = \rho^{1,1} : \KK^{1,1}(P) \to \KK^{1,1}(P)$ the birational rowmotion map on $P$.
We derive the reciprocity for $P$ from that for 
the rectangle poset $\tilde{P} = \{ (i,j) \in \Int^2 : 0 \le i, j \le r \}$.
We denote by $\tilde{\rho} : \KK^{1,1}(\tilde{P}) \to \KK^{1,1}(\tilde{P})$ 
the birational rowmotion map on $\tilde{P}$ with $A=B=1$.
The following lemma is a consequence of \cite[Lemma~59 (c)]{GR2} 
and Lemma~\ref{lem:A=B=1} (with $A = 1/2$ and $B=2$).

\begin{lemma}
\label{lem:doubling}
For $F \in \KK^{1,1}(P)$, we define $\tilde{F} \in \KK^{1,1}(\tilde{P})$ by
$$
\tilde{F}(i,j)
 =
\begin{cases}
 F(i,j) &\text{if $i \le j$,} \\
 F(j,i) &\text{if $i > j$.}
\end{cases}
$$
Then we have
$$
\left( \rho^k F \right) (i,j)
 =
\left( \tilde{\rho}^k \tilde{F} \right) (i,j)
\times
\begin{cases}
1/2 &\text{if $1 \le k \le i+j$,} \\
1 &\text{if $k=i+j+1$,} \\
2 &\text{if $i+j+2 \le k \le 2r+1$,} \\
1 &\text{if $k=2r+2$}
\end{cases}
$$
for $1 \le k \le 2r+2$ and $(i,j) \in P$.
\end{lemma}

By using this lemma and the reciprocity for the rectangle poset $\tilde{P}$, we have
$$
\left( \rho^{i+j+1} F \right)(i,j)
 =
\left( \tilde{\rho}^{i+j+1} \tilde{F} \right)(i,j)
 =
\frac{ 1 }{ \tilde{F} (r-i,r-j) }
 =
\frac{ 1 }{ F(r-j,r-i) }.
$$
This is the desired identity for a shifted staircase poset.

\subsection{%
Double-tailed diamond posets
}

In this subsection, we prove the reciprocity for double-tailed diamond posets.
Let $P = P_{D_n,\varpi_1}$ be the minuscule poset associated to the minuscule weight 
$\lambda = \varpi_1$ of the Lie algebra of type $D_n$.
We label elements of $P$ by
\begin{gather*}
v_i = e_1 + e_{i+1} \quad (1 \le i \le n-2),
\\
v_{n-1}^+ = e_1 + e_n,
\quad
v_{n-1}^- = e_1 - e_n,
\\
v_i = e_ 1 - e_{2n-1-i} \quad (n \le i \le 2n-3).
\end{gather*}
Note that $v_1$ is the maximum element and $v_{2n-3}$ is the minimum element.

Fix an initial state $X \in \KK^{1,1}(P)$.
We regard $X(v)$ ($v \in P$) as indeterminates and define $Z \in \KK^{1,1}(P)$ by (\ref{eq:X2Z}).
We write
\begin{gather*}
x_i = X(v_i) \quad(1 \le i \le 2n-3, \, i \neq n-1),
\quad
x_{n-1}^\pm = X(v_{n-1}^\pm),
\\
z_i = Z(v_i) \quad(1 \le i \le 2n-3, \, i \neq n-1),
\quad
z_{n-1}^\pm = Z(v_{n-1}^\pm).
\end{gather*}
Then we have
$$
z_i
 = 
\begin{cases}
 \dfrac{ x_i }{ x_{i+1} } &\text{if $i \neq n-1, 2n-3$,} \\
 \dfrac{ x_{n-2} }{ x_{n-1}^+ + x_{n-1}^- } &\text{if $i = n-1$,} \\
 x_{2n-3} &\text{if $i=2n-3$,}
\end{cases}
\quad
z_{n-1}^\pm = \frac{ x_{n-1}^\pm }{ x_n }.
$$
For positive integers $i$ and $l$ satisfying $1 \le i \le 2n-3$ and $i+l-1 \le 2n-3$, 
we define monomials $C(i;l)$ and $C^\pm(i;l)$ as follows:
\begin{enumerate}
\item[(i)]
If $1 \le i \le n-2$ and $i+l-1 \le n-2$, then we put
$$
C(i;l) = z_i z_{i+1} \cdots z_{i+l-1}.
$$
\item[(ii)]
If $1 \le i \le n-1$ and $n-1 \le i+l-1 \le 2n-3$, then we put
$$
C^\pm(i;l) = z_i z_{i+1} \cdots z_{n-2} z_{n-1}^\pm z_n \cdots z_{i+l-1}.
$$
\item[(iii)]
If $n+2 \le i \le 2n-3$, then we put
$$
C(i;l) = z_i z_{i+1} \cdots z_{i+l-1}.
$$
\end{enumerate}
Then the original indeterminates $X(v)$ can be expressed in terms of $Z(v)$ as follows:

\begin{lemma}
\label{lem:diamond}
The values $X(v)$ ($v \in P$) are expressed in terms of $C(i;l)$ and $C^\pm(i;l)$ as follows:
$$
\begin{cases}
X(v_i) = C^+(i;2n-i-2) + C^-(i;2n-i-2) &\text{if $1 \le i \le n-2$,} \\
X(v_{n-1}^\pm) = C^\pm(n-1;n-1) &\text{if $i=n-1$,} \\
X(v_i) = C(i;2n-i-2) &\text{if $n \le i \le 2n-3$.}
\end{cases}
$$
\end{lemma}

Recall that $P$ is a graded poset with rank function $\rank$ given by 
$\rank(v_i) = 2n-i-2$ ($1 \le i \le 2n-3$, $i \neq n-1$) and $\rank(v_{n-1}^\pm) = n-1$.
Then it is straightforward to prove the following explicit formulas 
by using the induction on $k$ and $i$.
(We omit the proof.)

\begin{prop}
\label{prop:diamond}
Let $v \in P$ and $k$ a positive integer.
If $1 \le k \le \rank(v)$, then the value $\left( \rho^k X \right)(v)$ of iterations of birational rowmotion 
is expressed in terms of $C(i;l)$ and $C^\pm(i;l)$ as follows:
\begin{enumerate}
\item[(a)]
If $v = v_i$ with $1 \le i \le n-2$, we have
$$
\left( \rho^k X \right) (v_i)
 =
\begin{cases}
 \dfrac{ 1 }{ C(k ; i) } &\text{if $1 \le k \le n-i-1$,} \\
 \dfrac{ 1 }{ C^+(k ; i) } + \dfrac{ 1 }{ C^-(k+1 ; i) } &\text{if $n-i \le k \le n-1$,} \\
 \dfrac{ 1 }{ C(k ; i) } &\text{if $n \le k \le 2n-i-2$.}
\end{cases}
$$
\item[(b)]
If $v = v_{n-1}^\pm$, we have
$$
\left( \rho^k X \right) (v_{n-1}^\pm)
 =
\frac{ 1 }
     { C^{\ep (-1)^{k-1}}(k ; n-1) }.
$$
\item[(c)]
If $v = v_i$ with $n \le i \le 2n-3$, we have
$$
\left( \rho^k X \right) (v_i)
 =
\frac{ 1 }
     { C^+(k ; i) + C^-(k ; i) }.
$$
\end{enumerate}
\end{prop}

Since the involution $\iota : P \to P$ is given by
$$
\iota(v_i) = v_{2n-i-2}
\quad(1 \le i \le 2n-3, \, i \neq n-1),
\quad
\iota(v_{n-1}^\ep) = v_{n-1}^{\ep (-1)^n},
$$
we obtain the desired reciprocity by comparing formulas in Lemma~\ref{lem:diamond} 
and Proposition~\ref{prop:diamond}.
This completes the proof of Theorem~\ref{thm:main1} (b) for all minuscule posets.

\subsection{%
Reciprocity for birational Coxeter-motion
}

We have the following conjectural reciprocity for a particular birational Coxeter-motion map.

\begin{conjecture}
\label{conj:half-period}
Let $P$ be a minuscule poset.
We decompose the simple root system $\Pi$ into a disjoin union of two subsets $\Pi_1$ and $\Pi_2$ 
such that any roots in $\Delta_i$ are pairwise orthogonal for each $i$. 
We define $\gamma_1$ and $\gamma_2$ by
$$
\gamma_1 = \prod_{\alpha \in \Pi_1} \sigma^{A,B}_\alpha,
\quad
\gamma_2 = \prod_{\beta \in \Pi_2} \sigma^{A,B}_\beta,
$$
and put
$$
\delta
 =
\underbrace{\gamma_1 \gamma_2 \gamma_1 \gamma_2 \gamma_1 \cdots}_{\text{$h$ factors}},
$$
where $h$ is the Coxeter number.
Then we conjecture that
\begin{equation}
\label{eq:half-period}
( \delta F )(v)
 =
\frac{ AB }
     { F( \iota v) }
\end{equation}
for any $F \in \KK^{A,B}(P)$ and $v \in P$.
\end{conjecture}

The periodicity of birational Coxeter-motion maps is a consequence of this conjecture.
In fact, $\gamma = \gamma_1 \gamma_2$ is a Coxeter-motion map and
$$
\gamma^h
 = 
\begin{cases}
 \delta^2 &\text{if $h$ is even,} \\
 \delta_{1,2} \delta_{2,1} &\text{if $n$ is odd,}
\end{cases}
$$
where $\delta_{1,2} = \gamma_1 \gamma_2 \gamma_1 \gamma_2 \gamma_1 \cdots \gamma_1$ 
and $\delta_{2,1} = \gamma_2 \gamma_1 \gamma_1 \gamma_2 \gamma_1 \cdots \gamma_2$. 
If $h$ is even, then we have
$$
\left( \gamma^h F \right) (v)
 =
\left( \delta^2 F \right) (v)
 =
\frac{ AB }
     { \left( \delta F \right) (\iota v) }
 =
\frac{ AB }
     { AB / F( \iota^2 v ) }
 =
F(v).
$$
If $h$ is odd, we can derive $\left( \gamma^h F \right)(v) = F(v)$ from (\ref{eq:half-period}) 
in a similar manner.
\section{%
File homomesy
}

This section is devoted to the proof of the file homomesy phenomenon 
(Theorem~\ref{thm:main1} (c) and Theorem~\ref{thm:main2} (b)).

\subsection{%
Local properties
}

First we investigate local properties of birational rowmotion 
and Coxeter-motion around a given file.

Let $P$ be a minuscule poset with coloring $c : P \to \Pi$.
We regard the Hasse diagram of the poset $\hat{P} = P \sqcup \{ \hat{1}, \hat{0} \}$ as a directed graph, 
where a directed edge $u \to v$ corresponds to the covering relation $u \lessdot v$.
For $\alpha \in \Pi$, let $\hat{N}^\alpha$ be the neighborhood of 
$P^\alpha = \{ x \in P : c(x) = \alpha \}$ given by
$$
\hat{N}^\alpha
 =
\{ x \in \hat{P} : \text{there is an element $y \in P^\alpha$ such that $x \lessdot y$ or $x \gtrdot y$} \}.
$$
We define $G^\alpha$ to be the bipartite directed subgraph of the Hasse diagram of $\hat{P}$ 
with black vertex set $P^\alpha$ and white vertex set $\hat{N}^\alpha$.
It follows from Proposition~\ref{prop:minuscule2} (c) that
$$
\hat{N}^\alpha
 =
\bigsqcup_{\beta \sim \alpha} P^\beta
 \sqcup
\begin{cases}
 \{ \hat{1}, \hat{0} \} &\text{if $\alpha = \alpha_{\max} = \alpha_{\min}$,} \\
 \{ \hat{1} \} &\text{if $\alpha = \alpha_{\max} \neq \alpha_{\min}$,} \\
 \{ \hat{0} \} &\text{if $\alpha = \alpha_{\min} \neq \alpha_{\max}$,} \\
 \emptyset &\text{otherwise,}
\end{cases}
$$
where $\beta$ runs over all simple roots adjacent to $\alpha$ in the Dynkin diagram, 
and $\alpha_{\max}$ (resp. $\alpha_{\min}$) is the color of the maximum (resp. minimum) 
element of $P$.

To describe the graph structure of $G^\alpha$, 
we introduce two series of posets $G_m$ and $H_m$.
For a positive integer $m$, let $G_m$ be the poset consisting of $3m$ elements 
$x_1, \cdots, x_m, y_1, \cdots, y_{m-1}, z_1, \allowbreak \cdots, z_{m-1}, u, v$ with covering relations
$$
u \lessdot x_1,
\quad
x_i \lessdot y_i \lessdot x_{i+1},
\quad
x_i \lessdot z_i \lessdot x_{i+1},
\quad
x_m \lessdot v.
$$
Note that $G_1$ is the three-element chain.
And, for an integer $m \ge 2$, let $H_m$ be the $(2m+1)$-element chain 
$$
u \lessdot x_1 \lessdot y_1 \lessdot x_2 \lessdot y_2 \lessdot \cdots \lessdot y_{m-1} \lessdot x_m \lessdot v
$$
We regard the Hasse diagrams of $G_m$ and $H_m$ as bipartite directed graphs with black vertices $x_1, \dots, x_m$. 
For example, the Hasse diagrams of $G_4$ and $H_4$ are shown in 
Figures~\ref{fig:G} and \ref{fig:H} respectively.
\begin{figure}[ht]
\setlength{\unitlength}{1.5pt}
\centering
\begin{minipage}{0.3\hsize}
\centering
\begin{picture}(50,90)
\put(25,5){\circle{3}}
\put(25,15){\circle*{3}}
\put(15,25){\circle{3}}
\put(35,25){\circle{3}}
\put(25,35){\circle*{3}}
\put(15,45){\circle{3}}
\put(35,45){\circle{3}}
\put(25,55){\circle*{3}}
\put(15,65){\circle{3}}
\put(35,65){\circle{3}}
\put(25,75){\circle*{3}}
\put(25,85){\circle{3}}
\put(15,65){\circle{3}}
\put(15,65){\circle{3}}
\put(25,6.5){\vector(0,1){7}}
\put(24,16){\vector(-1,1){8}}
\put(26,16){\vector(1,1){8}}
\put(16,26){\vector(1,1){8}}
\put(34,26){\vector(-1,1){8}}
\put(24,36){\vector(-1,1){8}}
\put(26,36){\vector(1,1){8}}
\put(16,46){\vector(1,1){8}}
\put(34,46){\vector(-1,1){8}}
\put(24,56){\vector(-1,1){8}}
\put(26,56){\vector(1,1){8}}
\put(16,66){\vector(1,1){8}}
\put(34,66){\vector(-1,1){8}}
\put(25,76.5){\vector(0,1){7}}
\put(25,80){\makebox(10,10){$v$}}
\put(25,70){\makebox(15,10){$x_4$}}
\put(25,50){\makebox(15,10){$x_3$}}
\put(25,30){\makebox(15,10){$x_2$}}
\put(25,10){\makebox(15,10){$x_1$}}
\put(25,0){\makebox(10,10){$u$}}
\put(0,60){\makebox(15,10){$y_3$}}
\put(0,40){\makebox(15,10){$y_2$}}
\put(0,20){\makebox(15,10){$y_1$}}
\put(35,60){\makebox(15,10){$z_3$}}
\put(35,40){\makebox(15,10){$z_2$}}
\put(35,20){\makebox(15,10){$z_1$}}
\end{picture}
\caption{$G_4$}
\label{fig:G}
\end{minipage}
\begin{minipage}{0.3\hsize}
\centering
\begin{picture}(50,90)
\put(25,5){\circle{3}}
\put(25,15){\circle*{3}}
\put(15,25){\circle{3}}
\put(25,35){\circle*{3}}
\put(15,45){\circle{3}}
\put(25,55){\circle*{3}}
\put(15,65){\circle{3}}
\put(25,75){\circle*{3}}
\put(25,85){\circle{3}}
\put(15,65){\circle{3}}
\put(15,65){\circle{3}}
\put(25,6.5){\vector(0,1){7}}
\put(24,16){\vector(-1,1){8}}
\put(16,26){\vector(1,1){8}}
\put(24,36){\vector(-1,1){8}}
\put(16,46){\vector(1,1){8}}
\put(24,56){\vector(-1,1){8}}
\put(16,66){\vector(1,1){8}}
\put(25,76.5){\vector(0,1){7}}
\put(25,80){\makebox(10,10){$v$}}
\put(25,70){\makebox(15,10){$x_4$}}
\put(25,50){\makebox(15,10){$x_3$}}
\put(25,30){\makebox(15,10){$x_2$}}
\put(25,10){\makebox(15,10){$x_1$}}
\put(25,0){\makebox(10,10){$u$}}
\put(0,60){\makebox(15,10){$y_3$}}
\put(0,40){\makebox(15,10){$y_2$}}
\put(0,20){\makebox(15,10){$y_1$}}
\end{picture}
\caption{$H_4$}
\label{fig:H}
\end{minipage}
\end{figure}

\begin{lemma}
\label{lem:decomp}
Each bipartite directed graph $G^\alpha$ is decomposed into a disjoint union of graphs of the form $G_m$ or $H_m$ 
as follows:
\begin{itemize}
\item
If $P = P_{A_n,\varpi_r}$, then
$$
G^{\alpha_i}
 \cong
\begin{cases}
 G_i &\text{if $1 \le i \le r$,} \\
 G_r &\text{if $r \le i \le s$,} \\
 G_{n-i+1} &\text{if $s \le i \le l$,}
\end{cases}
$$
where $r+s = n+1$ and $r \le s$.
\item
If $P = P_{B_n,\varpi_n}$, then
$$
G^{\alpha_i}
 \cong
\begin{cases}
 G_i &\text{if $1 \le i \le n-1$,} \\
 H_l &\text{if $i = n$.}
\end{cases}
$$
\item
If $P = P_{C_n,\varpi_1}$, then
$$
G^{\alpha_i}
 \cong
\begin{cases}
 G_1 \sqcup G_1 &\text{if $1 \le i \le n-2$,} \\
 H_2 &\text{if $i = n-1$,} \\
 G_1 &\text{if $i = n$.}
\end{cases}
$$
\item
If $P = P_{D_n,\varpi_1}$, then
$$
G^{\alpha_i}
 \cong
\begin{cases}
 G_1 \sqcup G_1 &\text{if $1 \le i \le n-3$,} \\
 G_2 &\text{if $i = n-2$,} \\
 G_1 &\text{if $i = n-1$, $n$.}
\end{cases}
$$
\item
If $P = P_{D_n,\varpi_n}$, then
$$
G^{\alpha_i}
 \cong
\begin{cases}
 G_l &\text{if $1 \le i \le n-2$,} \\
 \left( G_1 \right)^{\sqcup \lfloor (n-1)/2 \rfloor} &\text{if $i = n-1$,} \\
 \left( G_1 \right)^{\sqcup \lfloor n/2 \rfloor} &\text{if $i = n$,} \\
\end{cases}
$$
where $G_1^{\sqcup m}$ is the disjoint union of $m$ copies of $G_1$, 
and $\lfloor x \rfloor$ stands for the largest integer not exceeding $x$.
\item
If $P = P_{E_6,\varpi_6}$, then
$$
G^{\alpha_i}
 \cong
\begin{cases}
 G_1 \sqcup G_1 &\text{if $i=1$, $2$, $6$,} \\
 G_1 \sqcup G_2 &\text{if $i=3$, $5$,} \\
 G_4 &\text{if $i=4$.}
\end{cases}
$$
\item
If $P = P_{E_7, \varpi_7}$, then
$$
G^{\alpha_i}
 \cong
\begin{cases}
 G_1 \sqcup G_1 &\text{if $i=1$,} \\
 G_1 \sqcup G_1 \sqcup G_1 &\text{if $i=2$,} \\
 G_2 \sqcup G_2 &\text{if $i=3$,} \\
 G_6 &\text{if $i=4$,} \\
 G_1 \sqcup G_3 \sqcup G_1 &\text{if $i=5$,} \\
 G_1 \sqcup G_2 \sqcup G_1 &\text{if $i=6$,} \\
 G_1 \sqcup G_1 \sqcup G_1 &\text{if $i=7$.}
\end{cases}
$$
\end{itemize}
\end{lemma}

The following relations are a key to the proof of the file homomesy phenomenon.

\begin{lemma}
\label{lem:local1}
Let $\rho = \rho^{A,B}$ be the birational rowmotion map and $\alpha$ a simple root.
\begin{enumerate}
\item[(a)]
If $G_m$ appears as a connected component of $G^\alpha$, then we have
\begin{equation}
\label{eq:local1}
\prod_{i=1}^m (\rho^{i-1} F)(x_i) \cdot \prod_{i=1}^m (\rho^i F)(x_i)
=
F(u)
\cdot
\prod_{l=1}^{m-1} (\rho^i F)(y_i)
\cdot
\prod_{l=1}^{m-1} (\rho^i F)(z_i)
\cdot
(\rho^m F)(v).
\end{equation}
\item[(b)]
If $H_m$ appears as a connected component of $G^\alpha$, then we have
\begin{equation}
\label{eq:local3}
\prod_{i=1}^m (\rho^{i-1} F)(x_i) \cdot \prod_{i=1}^m (\rho^i F)(x_i)
 =
F(u)
\cdot
\prod_{l=1}^{m-1} (\rho^i F)(y_i)^2
\cdot
(\rho^m F)(v).
\end{equation}
\end{enumerate}
\end{lemma}

\begin{demo}{Proof}
(a)
It follows from (\ref{eq:Brow_inductive}) that
\begin{gather*}
F(x_m) \cdot (\rho F)(x_m)
 =
(\rho F)(v) \cdot ( F(y_{m-1}) + F(z_{m-1}) ),
\\
F(x_i) \cdot (\rho F)(x_i)
 =
\frac{ (\rho F) (y_i) \cdot (\rho F) (z_i) \cdot 
       ( F(y_{i-1}) + F(z_{i-1}) ) }
     { (\rho F)(y_i) + (\rho F)(z_i) }
\quad(2 \le i \le m-1),
\\
F(x_1) \cdot (\rho F)(x_1)
 =
\frac{ F(u) \cdot (\rho F)(y_1) \cdot (\rho F)(z_1) }
     { (\rho F)(y_1) + (\rho F)(z_1) }.
\end{gather*}
By replacing $F$ with $\rho^{m-1} F$ (resp. $\rho^{i-1} F$) in the first (resp. second) equation, 
and then by multiplying the resulting equations together, we obtain (\ref{eq:local1}).

(b) can be checked by a similar computation.
\qed
\end{demo}

\begin{lemma}
\label{lem:local2}
Let $\alpha$ be a simple root and 
and $\sigma_\alpha = \prod_{v \in P^\alpha} \tau^{A,B}_v$ 
the product of birational toggles over $P^\alpha$.
\begin{enumerate}
\item[(a)]
If $G_m$ appears as a connected component of $G^\alpha$, then we have
\begin{equation}
\label{eq:local2}
\prod_{i=1}^m F(x_i) \cdot \prod_{i=1}^m (\sigma_\alpha F)(x_i)
=
F(u)
\cdot
\prod_{l=1}^{m-1} F(y_i)
\cdot
\prod_{l=1}^{m-1} F(z_i)
\cdot
F(v).
\end{equation}
\item[(b)]
If $H_m$ appears as a connected component of $G^\alpha$, then we have
\begin{equation}
\label{eq:local4}
\prod_{i=1}^m F(x_i) \cdot \prod_{i=1}^m (\sigma_\alpha F)(x_i)
 =
F(u)
\cdot
\prod_{l=1}^{m-1} F(y_i)^2
\cdot
F(v).
\end{equation}
\end{enumerate}
\end{lemma}

\begin{demo}{Proof}
(a)
By the definition (\ref{eq:Btoggle}), we have
\begin{gather*}
F(x_m) \cdot (\sigma_\alpha F)(x_m)
 =
F(v) \cdot ( F(y_{m-1}) + F(z_{m-1}) ),
\\
F(x_i) \cdot (\sigma_\alpha F)(x_i)
 =
\frac{ F(y_i) \cdot F(z_i) \cdot ( F(y_{i-1}) + F(z_{i-1}) ) }
     { F(y_i) + F(z_i) }
\quad(2 \le i \le m-1),
\\
F(x_1) \cdot (\sigma_\alpha F)(x_1)
 =
\frac{ F(y_1) \cdot F(z_1) \cdot F(u) }
     { F(y_1) + F(z_1) }.
\end{gather*}
Multiplying them together, we obtain (\ref{eq:local2}).

(b) can be checked by a similar computation.
\qed
\end{demo}

\subsection{%
File homomesy for birational rowmotion
}

In this subsection, we prove the file homomesy phenomenon for birational rowmotion 
(Theorem~\ref{thm:main1} (c)).

The following properties of Coxeter elements will be useful 
in the proof of Theorem~\ref{thm:main1} (c) and Theorem~\ref{thm:main2}(b);
the proof of the latter will be given in the next subsection. 
A \emph{Coxeter element} in a Weyl group $W = \langle s_\alpha : \alpha \in \Pi \rangle$ 
is a product of all simple reflections $s_\alpha$ in any order.
Then it is known that all Coxeter elements are conjugate.
By definition, the Coxeter number is the order of any Coxeter element.

\begin{lemma}
\label{lem:Coxeter}
Let $c$ be a Coxeter element and $h$ the Coxeter number.
Then we have
\begin{enumerate}
\item[(a)]
If $\mu \in \mathfrak{h}^*$ satisfies $c \mu = \mu$, then $\mu = 0$.
\item[(b)]
As a linear transformation on $\mathfrak{h}^*$,  we have
\begin{equation}
\label{eq:Coxeter1}
\sum_{k=0}^{h-1} c^k = 0
\end{equation}
\item[(c)]
Let $\alpha \in \Pi$ be a simple root and $\varpi$ the corresponding fundamental weight.
If $c = s_{\alpha_1} \cdots s_{\alpha_n}$ is a Coxeter element with $\Pi = \{ \alpha_1, \dots, \alpha_n \}$ 
and $\beta = s_{\alpha_1} \cdots s_{\alpha_{k-1}} \alpha_k$, where $\alpha = \alpha_k$, then we have
\begin{gather}
\label{eq:Coxeter2}
c \varpi = \varpi - \beta,
\\
\label{eq:Coxeter3}
\sum_{k=1}^{h-1} \sum_{i=0}^{k-1} c^i ( \beta )
 =
h \varpi.
\end{gather}
\end{enumerate}
\end{lemma}

\begin{demo}{Proof}
(a)
See \cite[V, \S6, n${}^\circ$2]{B1}.

(b) follows from $c^h = 1$ and (a).

(c)
Since $s_\gamma \varpi = \varpi - \langle \gamma^\vee, \varpi \rangle \gamma
 = \varpi - \delta_{\alpha,\gamma} \alpha$ for $\gamma \in \Pi$, 
we have 
$c \varpi
 =
\varpi - s_{\alpha_1} \cdots s_{\alpha_{k-1}} \alpha_k
 =
\varpi - \beta$.
Hence we see that
$$
c^k \varpi = \varpi - \sum_{i=0}^{k-1} c^i \beta.
$$
By using (\ref{eq:Coxeter1}), we obtain
$$
0
 =
\sum_{k=0}^{h-1} c^k \varpi
 =
h \varpi - \sum_{k=1}^{h-1} \sum_{i=0}^{k-1} c^i \beta,
$$
from which (\ref{eq:Coxeter3}) follows.
\qed
\end{demo}

In order to prove Theorem~\ref{thm:main1} (c), we consider
\begin{equation}
\label{eq:Phi'}
\Phi'_\alpha(F)
 =
\prod_{v \in P^\alpha}
 \left( \rho^{(\rank(v) - \rank(v^\alpha_0))/2} F \right)(v),
\end{equation}
instead of $\Phi_\alpha(F) = \prod_{v \in P^\alpha} F(v)$.
Here $v^\alpha_0$ is the minimum element of $P^\alpha$.
Note that $P^\alpha$ is a chain and $\rank(v) - \rank(v^\alpha_0)$ is an even integer 
(see Proposition~\ref{prop:minuscule2} (d) and (e)).
Since $\rho$ has a finite order $h$, we have
\begin{equation}
\label{eq:Phi=Phi'}
\prod_{k=0}^{h-1} \Phi_\alpha( \rho^k F )
 =
\prod_{k=0}^{h-1} \Phi'_\alpha( \rho^k F ).
\end{equation}

\begin{remark}
It is worth mentioning that $\Phi'_\alpha( \rho^k X)$ are Laurent monomials in the variables $Z(v)$ 
defined by (\ref{eq:X2Z}).
In a forthcoming paper \cite{O}, 
we will give explicit formulas for $\Phi'_\alpha( \rho^k X)$ in classical types.
\end{remark}

\begin{prop}
\label{prop:rel_Phi'}
For $\alpha \in \Pi$ and $F \in \KK^{A,B}(P)$, we have
\begin{equation}
\label{eq:rel_Phi'}
\Phi'_\alpha(F) \cdot \Phi'_\alpha ( \rho F )
 =
A^{\delta_{\alpha,\alpha_{\max}}} B^{\delta_{\alpha,\alpha_{\min}}}
\prod_{\beta \sim \alpha} \Phi'_\beta ( \rho^{m_{\alpha,\beta}} F)^{-\langle \beta, \alpha^\vee \rangle},
\end{equation}
where $\beta$ runs over all simple roots adjacent to $\alpha$ in the Dynkin diagram and 
$$
m_{\alpha,\beta}
 =
\begin{cases}
 1 &\text{if $v^\beta_0 > v^\alpha_0$,} \\
 0 &\text{if $v^\beta_0 < v^\alpha_0$.}
\end{cases}
$$
\end{prop}

\begin{demo}{Proof}
We explain the proof in the case where $\mathfrak{g}$ is of type $E_7$, $\lambda = \varpi_7$ 
and $\alpha = \alpha_5$. (The other cases can be proved in a similar way.)
We label elements of $P^\alpha$ as $v^\alpha_0, v^\alpha_1, v^\alpha_2, \dots$ from bottom to top.
By the definition (\ref{eq:Phi'}), we have
\begin{align*}
\Phi'_{\alpha_4}(F)
 &=
F(v^{\alpha_4}_0) \cdot (\rho F)(v^{\alpha_4}_1) \cdot
(\rho^2 F)(v^{\alpha_4}_2) \cdot (\rho^3 F)(v^{\alpha_4}_3) \cdot
(\rho^4 F)(v^{\alpha_4}_4) \cdot (\rho^5 F)(v^{\alpha_4}_5),
\\
\Phi'_{\alpha_5}(F)
 &=
F(v^{\alpha_5}_0) \cdot (\rho^2 F)(v^{\alpha_5}_1) \cdot
(\rho^3 F)(v^{\alpha_5}_2) \cdot (\rho^4 F)(v^{\alpha_5}_3) \cdot
(\rho^6 F)(v^{\alpha_5}_4),
\\
\Phi'_{\alpha_6}(F)
 &=
F(v^{\alpha_6}_0) \cdot (\rho^3 F)(v^{\alpha_6}_1) \cdot 
(\rho^4 F)(v^{\alpha_6}_2) \cdot (\rho^7 F)(v^{\alpha_6}_3).
\end{align*}
The subgraph $G^{\alpha_5}$ has three connected components 
\begin{gather*}
\{ v^{\alpha_4}_0, v^{\alpha_5}_0, v^{\alpha_6}_0 \} \cong G_2,
\\
\{ v^{\alpha_4}_1, v^{\alpha_4}_2, v^{\alpha_4}_3, v^{\alpha_4}_4, 
v^{\alpha_5}_1, v^{\alpha_5}_2, v^{\alpha_5}_3, v^{\alpha_6}_1, v^{\alpha_6}_2 \}
 \cong G_3,
\\ 
\{ v^{\alpha_4}_5, v^{\alpha_5}_4, v^{\alpha_6}_3 \} \cong G_1.
\end{gather*} 
By applying (\ref{eq:local1}) to three connected components of $G^{\alpha_5}$, 
we obtain
\begin{align*}
&
F(v^{\alpha_5}_0) \cdot (\rho F)(v^{\alpha_5}_0)
 =
F(v^{\alpha_6}_0) \cdot (\rho F)(v^{\alpha_4}_0),
\\
&
F(v^{\alpha_5}_1) \cdot (\rho F)(v^{\alpha_5}_2) \cdot (\rho^2 F)(v^{\alpha_5}_3) \cdot 
(\rho F)(v^{\alpha_5}_1) \cdot (\rho^2 F)(v^{\alpha_5}_2) \cdot (\rho^3 F)(v^{\alpha_5}_3)
\\
&\quad
 =
F(v^{\alpha_4}_1) \cdot 
(\rho F)(v^{\alpha_6}_1) \cdot (\rho F)(v^{\alpha_4}_2) \cdot
(\rho^2 F)(v^{\alpha_6}_2) \cdot (\rho^2 F)(v^{\alpha_4}_3) \cdot
(\rho^3 F)(v^{\alpha_4}_4)
\\
&
F(v^{\alpha_5}_4) \cdot (\rho F)(v^{\alpha_5}_4)
 =
F(v^{\alpha_4}_5) \cdot (\rho F)(v^{\alpha_6}_3).
\end{align*}
By replacing $F$ with $\rho^2 F$ (resp. $\rho^6 F$) in the second (resp. third) equation, 
and then by multiplying three resulting equations together, we have
$$
\Phi'_{\alpha_5}(F) \cdot \Phi'_{\alpha_5}(\rho F)
 =
\Phi'_{\alpha_6}(F) \cdot \Phi'_{\alpha_4}(\rho F).
$$
Since $v^{\alpha_6}_0 < v^{\alpha_5}_0 < v^{\alpha_4}_0$ (see Figure~\ref{fig:e7}), 
we obtain (\ref{eq:rel_Phi'}) in this case.
\qed
\end{demo}

\begin{corollary}
\label{cor:rel_Phi}
For a simple root $\beta \in \Pi$, we put
$$
\tilde{\Phi}_\beta(F)
 = 
\prod_{k=0}^{h-1} \Phi_\beta(\rho^k F).
$$
Then we have for fixed $\alpha \in \Pi$,
\begin{equation}
\label{eq:rel_Phi}
\prod_{\beta \in \Pi} 
 \tilde{\Phi}_\beta(F)^{\langle \beta, \alpha^\vee \rangle}
 =
A^{\delta_{\alpha,\alpha_{\max}}} B^{\delta_{\alpha,\alpha_{\min}}}
\end{equation}
for any $F \in \KK^{A,B}(P)$.
\end{corollary}

\begin{demo}{Proof}
Since $\rho$ has a finite order $h$, Equation (\ref{eq:Phi=Phi'}) implies 
$\tilde{\Phi}_\beta(F) = \prod_{k=0}^{h-1} \Phi'_\beta(\rho^{k+m} F)$ for any integer $m$.
Hence (\ref{eq:rel_Phi}) follows from (\ref{eq:rel_Phi'}).
\qed
\end{demo}

Now we are ready to prove Theorem~\ref{thm:main1} (c).

\begin{demo}{Proof of Theorem~\ref{thm:main1} (c)}
We define an element $\tilde{\mu}(F) \in \mathfrak{h}^*$ for $F \in \KK^{A,B}(P)$ by putting
$$
\tilde{\mu}(F)
 =
\sum_{\alpha \in \Pi} \log \tilde{\Phi}_\alpha(F) \cdot \alpha.
$$
Note that, if $\varpi^\vee$ is the fundamental coweight corresponding to $\alpha$, 
then we have
$$
\log \tilde{\Phi}_\alpha(F) = \langle \varpi^\vee, \tilde{\mu}(F) \rangle.
$$
Since $\varpi_{\max} = - w_0 \lambda$ (resp. $\varpi_{\min} = \lambda$) 
is the fundamental weight corresponding to the color $\alpha_{\max}$ (resp. $\alpha_{\min}$) 
of the maximum (resp. minimum) element of $P$ (see Proposition~\ref{prop:minuscule2} (b)), 
it it enough to show 
\begin{equation}
\label{eq:mu}
\tilde{\mu}(F)
 =
h a \cdot \varpi_{\max} + h b \cdot \varpi_{\min},
\end{equation}
where $a = \log A$, $b = \log B$.

Since we have
$$
\sum_{\beta \in \Pi} \langle \beta, \alpha^\vee \rangle \log \tilde{\Phi}_\beta(F)
=
h a \delta_{\alpha,\alpha_{\max}}
+
h b \delta_{\alpha,\alpha_{\max}}
$$
by Corollary~\ref{cor:rel_Phi}, we see that for any $\alpha \in \Pi$ 
\begin{align*}
s_\alpha \tilde{\mu}(F)
 &=
\sum_{\beta \in \Pi}
 \log \tilde{\Phi}_\beta(F) \cdot ( \beta - \langle \beta, \alpha^\vee \rangle \alpha )
\\
 &=
\sum_{\beta \in \Pi} \log \tilde{\Phi}_\beta(F) \beta
-
\left( \sum_{\beta \in \Pi} \langle \beta, \alpha^\vee \rangle \log \tilde{\Phi}_\beta(F) \right) \alpha
\\
 &=
\tilde{\mu}(F)
 -
\left( \delta_{\alpha,\alpha_{\max}} h a + \delta_{\alpha,\alpha_{\min}} h b \right) \alpha.
\end{align*}
Let $c = s_{\alpha_1} \cdots s_{\alpha_n}$ be a Coxeter element 
and put
$$
\beta_{\max} = s_{\alpha_1} \cdots s_{\alpha_{k-1}} \alpha_k,
\quad 
\beta_{\min} = s_{\alpha_1} \cdots s_{\alpha_{m-1}} \alpha_m, 
$$
where $\alpha_k = \alpha_{\max}$, $\alpha_m = \alpha_{\min}$.
Then we have
$$
c \tilde{\mu}(F)
 =
\tilde{\mu}(F)
 -
\left( h a \cdot \beta_{\max} + h b \cdot \beta_{\min} \right).
$$
By substituting $\beta_{\max} = \varpi_{\max} - c \varpi_{\max}$ 
and $\beta_{\min} = \varpi_{\min} -c \varpi_{\min}$ (see (\ref{eq:Coxeter2})), 
we have
$$
c \left( \tilde{\mu}(F) - h a \cdot \varpi_{\max} - h b \cdot \varpi_{\max} \right)
 =
\tilde{\mu}(F) - h a \cdot \varpi_{\max} - h b \cdot \varpi_{\max}.
$$
Then it follows from Lemma~\ref{lem:Coxeter} (a) that
$$
\tilde{\mu}(F) - h a \cdot \varpi_{\max} - h b B \cdot \varpi_{\max} = 0.
$$
This completes the proof of (\ref{eq:mu}) and hence of Theorem~\ref{thm:main1} (c).
\qed
\end{demo}

\subsection{%
File homomesy for birational Coxeter-motion
}

In this subsection we prove Theorem~\ref{thm:main2} (b).
The following proposition is a consequence of Lemma~\ref{lem:decomp} 
and Equations (\ref{eq:local2}), (\ref{eq:local4}).

\begin{prop}
\label{prop:rel_Phi}
Let $\sigma_\alpha = \prod_{v \in P^\alpha} \tau_v : \KK^{A,B}(P) \to \KK^{A,B}(P)$ 
be the product of toggles over $P^\alpha$.
Then we have
\begin{enumerate}
\item[(a)]
For a simple root $\alpha$, we have
$$
\Phi_\alpha(F) \cdot \Phi_\alpha ( \sigma_\alpha F )
 =
A^{\delta_{\alpha,\alpha_{\max}}} B^{\delta_{\alpha,\alpha_{\min}}}
\prod_{\beta \sim \alpha} \Phi_\beta(F)^{-\langle \beta, \alpha^\vee \rangle}.
$$
\item[(b)]
For simple roots $\alpha \neq \beta$, we have $\Phi_\beta (\sigma_\alpha F) = \Phi_\beta(F)$.
\end{enumerate}
\end{prop}

By using this proposition, we can complete the proof of the file homomesy phenomenon 
for birational Coxeter-motion.

\begin{demo}{Proof of Theorem~\ref{thm:main2} (b)}
We define an element $\mu(F) \in \mathfrak{h}^*$ for $F \in \KK^{A,B}(P)$ by putting
$$
\mu(F) = \sum_{\beta \in \Pi} \log \Phi_\beta(F) \cdot \beta.
$$

First we prove 
\begin{equation}
\label{eq:s-mu}
\mu (\sigma_\alpha F)
 = 
s_\alpha \mu(F)
 +
\left(
 \delta_{\alpha,\alpha_{\max}} a + \delta_{\alpha,\alpha_{\min}} b
\right) \alpha
\end{equation}
where $a = \log A$ and $b = \log B$.
By using Proposition~\ref{prop:rel_Phi}, we have
\begin{align*}
\mu (\sigma_\alpha F)
 &=
\sum_{\beta \neq \alpha} \log \Phi_\beta (\sigma_\alpha F) \beta 
+
\log \Phi_\alpha (\sigma_\alpha F) \alpha
 \\
 &=
\sum_{\beta \neq \alpha} \log \Phi_\beta (F) \beta 
+
\left(
 \delta_{\alpha,\alpha_{\max}} a + \delta_{\alpha,\alpha_{\min}} b 
 - \sum_{\beta \neq \alpha} \langle \beta, \alpha^\vee \rangle \log \Phi_\beta(F)
 - \log \Phi_\alpha(F)
\right) \alpha
 \\
 &=
\sum_{\beta \neq \alpha} \log \Phi_\beta(F) (\beta - \langle \beta, \alpha^\vee \rangle \alpha)
 - \log \Phi_\alpha(F) \alpha
+
\left(
 \delta_{\alpha,\alpha_{\max}} a + \delta_{\alpha,\alpha_{\min}} b 
\right) \alpha
\\
 &=
\sum_{\beta \neq \alpha} \log \Phi_\beta(F) s_\alpha(\beta)
+ \log \Phi_\alpha(F) s_\alpha(\alpha)
+
\left(
 \delta_{\alpha,\alpha_{\max}} a + \delta_{\alpha,\alpha_{\min}} b 
\right) \alpha
\\
 &=
s_\alpha(\mu(F)) 
+
\left(
 \delta_{\alpha,\alpha_{\max}} a + \delta_{\alpha,\alpha_{\min}} b 
\right) \alpha.
\end{align*}

Suppose that $\gamma = \sigma_{\alpha_1} \cdots \sigma_{\alpha_n}$, 
and let $c = s_{\alpha_1} \cdots s_{\alpha_n}$ be the corresponding Coxeter element.
Then, by iteratively using (\ref{eq:s-mu}), we obtain
$$
\mu( \gamma F )
 =
c ( \mu(F) ) + a \cdot \beta_{\max} + b \cdot \beta_{\min},
$$
where $\beta_{\max}$ and $\beta_{\min}$ are defined by 
$\beta_{\max} = s_{\alpha_1} \cdots s_{\alpha_{k-1}} \alpha_k$, 
$\beta_{\min} = s_{\alpha_1} \cdots s_{\alpha_{m-1}} \alpha_m$ 
with $\alpha_k = \alpha_{\max}$ and $\alpha_m = \alpha_{\min}$.
Hence by induction on $k$ we see that
$$
\mu (\gamma^k F)
 = 
c^k(\mu(F))
 + 
a \sum_{i=0}^{k-1} c^i(\beta_{\max})
 +
b \sum_{i=0}^{k-1} c^i(\beta_{\min}).
$$
Therefore we have
$$
\sum_{k=0}^{h-1} \mu ( \gamma^k F )
 =
\sum_{k=0}^{h-1} c^k (\mu(F))
 +
a \sum_{k=1}^{h-1} \sum_{i=0}^{k-1} c^i(\beta_{\max})
 +
b \sum_{k=1}^{h-1} \sum_{i=0}^{k-1} c^i(\beta_{\min}).
$$
Now it follows from (\ref{eq:Coxeter1}) and (\ref{eq:Coxeter3}) that
$$
\sum_{k=0}^{h-1} \mu ( \gamma^k F )
 =
a h \cdot \varpi_{\max} + b h \cdot \varpi_{\min}.
$$
By the definition of $\mu(F)$, we have
$$
\sum_{\beta \in \Pi}
 \log \left( \prod_{k=0}^{h-1} \Phi_\beta( \gamma^k F) \right) \cdot \beta
 =
a h \cdot \varpi_{\max} + b h \cdot \varpi_{\min}.
$$
Then we can complete the proof by taking the pairing $\langle \quad, \varpi^\vee \rangle$.
\qed
\end{demo}



\end{document}